\documentclass{amsart}

\usepackage{amsthm,amsfonts,amsmath,amssymb,mathrsfs,verbatim,cite,graphicx,hyperref}
\pdfoutput=1

\newcommand{\la}{\lambda}
\newcommand{\Symm}{\mathfrak{S}}
\newcommand{\Poch}[3]{[#1]_{#2}^{(#3)}}
\newcommand{\R}{\mathbb R}
\newcommand{\Complex}{\mathbb C}
\newcommand{\Z}{\mathbb Z}
\DeclareMathOperator{\eup}{e}
\DeclareMathOperator{\height}{\textup{ht}}
\DeclareMathOperator{\alt}{\textup{alt}}
\newcommand{\dup}{\textup{d}}
\newcommand{\iup}{\textup{i}\hspace{1pt}}
\newcommand{\Int}{\int\limits}
\newcommand{\abs}[1]{\lvert#1\rvert}
\newcommand{\Abs}[1]{\bigl\lvert#1\bigr\rvert}
\renewcommand{\Re}{\textup{Re}}
\renewcommand{\Im}{\textup{Im}}
\newcommand*{\urlw}[1]{\href{http://www.#1}{\nolinkurl{www.#1}}}
\newcommand{\qbin}[2]{\genfrac{[}{]}{0pt}{}{#1}{#2}}
\newcommand{\g}{\mathfrak{g}}
\newcommand{\h}{\mathfrak{h}}
\newcommand{\ba}{\alpha}
\newcommand{\bil}[2]{(#1,#2)}
\newcommand{\gsl}{\mathfrak{sl}}
\renewcommand{\IJ}{I\hspace{-1pt}J}

\begin{document}

\title{The importance of the Selberg integral} 
\author{Peter J. Forrester and S. Ole Warnaar}\thanks{
This work has been supported by the Australian Research Council.
We gratefully acknowledge helpful correspondence on the mathematics
and history of the Selberg integral from 
George E.~Andrews, Enrico Bom\-bi\-eri,
Freeman J.~Dyson, Ron J.~Evans, Jon P.~Keating, Eric M.~Rains, Amitai Regev,
Vyacheslav P. Spiridonov and Richard P.~Stanley.}
\address{Department of Mathematics and Statistics, 
University of Melbourne, Victoria 3010, Australia}
\date{}

\begin{abstract} 
It has been remarked that a fair measure of the impact of Atle Selberg's
work is the number of mathematical terms which bear his name.
One of these is the Selberg integral, an $n$-dimensional
generalization of the Euler beta integral. We trace its sudden rise to
prominence, initiated by a question to Selberg from
Enrico Bombieri, more than thirty years after publication. 
In quick succession the Selberg integral was used to prove an 
outstanding conjecture in random matrix theory, and cases of the 
Macdonald conjectures. It further initiated the study of $q$-analogues, 
which in turn enriched the Macdonald conjectures.
We review these developments and proceed to exhibit the sustained 
prominence of the Selberg integral,
evidenced by its central role in random matrix theory, 
Calogero--Sutherland quantum many body systems,
Knizhnik--Zamolodchikov equations, and multivariable orthogonal
polynomial theory.
\end{abstract}

\subjclass[2000]{00-02,33-02}

\maketitle

\section{Discovery and reappearance}
\subsection*{1941 and 1944}
With the passing of Atle Selberg on August 6th 2007 at age 90, 
it is timely to reflect on his mathematical legacy. 
Indeed a number of brief articles highlighting some of his 
most influential mathematical discoveries were written shortly after the 
news of his death, see e.g.,~\cite{IAS07}. 
It is our aim to add to these tributes by giving a more comprehensive 
account of the mathematics, both pure and applied, related
to what now is referred to as the Selberg integral
\begin{align}\label{SelInt}
S_n(\alpha,\beta,\gamma) & :=
\int_0^1 \cdots \int_0^1 \,
 \prod_{i=1}^n t_i^{\alpha-1}(1-t_i)^{\beta-1}
\prod_{1\le i < j\le n} \abs{t_i - t_j}^{2\gamma}\,
\dup t_1\cdots\dup t_n \\
&\phantom{:}=
\prod_{j=0}^{n-1} \frac{\Gamma (\alpha+j\gamma)
\Gamma(\beta+j\gamma)\Gamma(1+(j+1)\gamma)}
{\Gamma(\alpha+\beta+(n+j-1)\gamma)\Gamma(1+\gamma)}.
\notag
\end{align}
The evaluation of this integral is valid for complex parameters 
$\alpha,\beta,\gamma$ such that
\begin{equation}\label{SelIntpv}
\Re(\alpha)>0,~\Re(\beta)>0,~\Re(\gamma)>
-\min\{1/n,\Re(\alpha)/(n-1),\Re(\beta)/(n-1)\},
\end{equation}
corresponding to the domain of convergence of the integral.

The proof of \eqref{SelInt} is the subject of Selberg's 
1944 paper ``Bemerkninger om et multipelt integral'' 
(Remarks on a multiple integral) \cite{Selberg44} --- the only one of
Selberg's works written in Norwegian --- published in 
\textit{Norsk Matematisk Tidsskrift}.
The latter has been compared \cite{Boas54} to the Scandinavian equivalent of 
the \textit{Mathematical Gazette}, with contents ranging from short 
research papers on subjects of general interest to discussions on 
teaching problems. Selberg himself remarks 
in his collected works \cite{Selberg89} that
\begin{quote}
This paper was published with some hesitation, and in Norwegian,
since I was rather doubtful that the results were new. The journal
is one which is read by mathematics-teachers in the gymnasium, and
the proof was written out in some detail so it should be understandable
to someone who knew a little about analytic functions and analytic 
continuation.
\end{quote}

Selberg's proof of \eqref{SelInt} proceeds by supposing $\gamma$ 
is a positive integer, and expanding
\[
\prod_{1\le i<j\le n} \abs{t_i-t_j}^{2\gamma} =
\sum_{0\le k_1,\dots,k_n \le 2(n-1)\gamma} c_{k_1,\dots,k_n}
t_1^{k_1} \cdots t_n^{k_n}.
\]
Substituting this expansion in the definition of $S_n(\alpha,\beta,\gamma)$
allows the resulting integrals to be evaluated by the Euler beta integral
\cite{Euler1730}
\begin{equation}\label{betaInt}
B(\alpha,\beta):=
\int_0^1 t^{\alpha-1} (1-t)^{\beta-1} \dup t =
\frac{\Gamma(\alpha)\Gamma(\beta)}{\Gamma(\alpha+\beta)},
\end{equation}
which itself is \eqref{SelInt} with $n=1$. 
The details, in English, of the proof from here on 
can be found in \cite{Mehta04,Fo02}, for example. 
Perhaps the most significant feature is the final step. 
It requires analytically continuing off the integers. 
Thus with \eqref{SelInt} established for $\gamma$ a positive integer, the
remaining task is to 
establish its validity for all 
complex $\gamma$ such that both sides are well defined.

For this purpose, after noting that both the left- and right-hand side
of \eqref{SelInt} are bounded analytic functions of $\gamma$ for 
$\Re(\gamma) \ge 1$ at least, Carlson's theorem can be used \cite{Carlson14}. 
The latter applies to functions $f(z)$ analytic for $\Re(z) \ge 0$ 
satisfying the bound $\abs{f(z)}=\text{O}(\eup^{\mu\abs{z}})$, $\mu < \pi$.
The theorem asserts that if furthermore $f(z)=0$ on the nonnegative
integers then, identically, $f(z)=0$. Note that the function
$f(z) = \sin \pi z$ shows that the bound $\mu < \pi$ is optimal. 
Selberg did not make direct use of  Carlson's theorem, but rather 
derived from first principles the same result in the case that $f(z)$ 
is bounded in the right half-plane, which is all that is required to
finalise the proof of \eqref{SelInt}.

Interestingly, although \cite{Selberg44} contains the first proof of 
\eqref{SelInt}, it is not the first time it appeared in print. 
This occurred three years earlier (albeit with the change of variables 
$t_i = s_i/(1+s_i)$ so that $s_i \in [0,\infty)$)
in Selberg's 1941 paper ``\"Uber einen Satz von A.~Gelfond'' 
(On a theorem of A.~Gelfond) \cite{Se41}.
Like \cite{Selberg44}, this earlier paper appeared in a Norwegian journal, 
this time \textit{Archiv for Mathematik og Naturvidenskab},
known for having Sophus Lie as one of its founders.
In a footnote Selberg remarks
\begin{quote}
Leider habe ich die Formel (11) [The Selberg integral] nirgends in der
Litteratur finden k\"onnen, ein Beweis hier zu bringen scheint aber nicht
angebracht, da die Arbeit sonst zu sehr anschwellen w\"urde; sollte sich aber
herausstellen, dass die Formel neu w\"are, beabsichtige ich sp\"ater ein Beweis
zu ver\"offentlichen.

\medskip

\noindent
[Unfortunately I have been unable to find formula (11) [The Selberg integral]
in the literature.  To present a proof here, however, seems inappropriate,
as it would make this paper significantly longer. If it turns out that the
formula is new, I intend to publish a proof at a later date.]
\end{quote}

Curiously, Selberg used his integral in \cite{Se41} to prove 
a result of some similarity to Carlson's theorem. As already noted,
the latter is itself an ingredient in Selberg's proof of \eqref{SelInt}.
Selberg's result relates to entire functions $f(z)$ of exponential type 
$\sigma(f)$, defined by
\[
\sigma(f) :=\mathop{\lim\sup}\limits_{r\to\infty} 
\frac{1}{r}
\log\Bigl(\:\mathop{\max}\limits_{\abs{z}=r}\abs{f(z)}\Bigr).
\]
A theorem of Hardy and P\'olya \cite{Polya00,Boas00} states that if 
$\sigma(f) < \log 2$ and $f$ takes integer values at the nonnegative 
integers, then $f(z)$ is polynomial. The transcendental function 
$f(z) = 2^z$ shows that this bound is optimal. 
A.~Gelfond \cite{Gelfond29} generalized this by proving that if 
$\sigma(f) < n \log (1+\exp(1/n-1))$ and $f$, together with its first 
$n-1$ derivatives, take integer values at the nonnegative integers,
then $f(z)$ is a polynomial. However, for $n>1$ this bound is
not optimal. By using his integral Selberg improved Gelfond's bound 
for $n>1$ to $\sigma(f) < \log m_n$, where $m_n$ is the minimum 
value of $\prod_{i=1}^n (1 + y_i)$ under the conditions $y_i >0$,
$y_1\cdots y_n =\eup^{1-n}$ and
$\prod_{1\le i<j\le n}\abs{y_i^{-1} - y_j^{-1}} = 1$.
This improves Gelfond's result since
\[
\prod_{i=1}^n (1+y_i) > 
(1 + \eup^{1/n-1} )^n.
\]

\subsection*{The 1950s to the late 1970s --- the Mehta integral}
For over thirty years the Selberg integral went essentially unnoticed. 
It was used only once --- in the special case $\alpha=\beta=1$, 
$\gamma=2$ --- in a study by S.~Karlin and L.S.~Shapley 
relating to the volume of a certain moment space, published in 1953 
\cite{KS53}.

Around ten years later there was again good reason to
make use of \eqref{SelInt}. Building upon the earlier work 
of E.P.~Wigner in the 1950s, F.J.~Dyson wrote a series of 
papers on the statistical theory of energy levels of complex systems. 
These papers ranged from the theory's foundations to its practical 
use in the analysis of experimental data. 
This last topic was addressed in part V of the series, in a work,
written jointly with M.L.~Mehta and published in 1963, which also summarizes 
both the status of the theory and open problems from that date.

A basic point is that random Hermitian matrices are used to model
the highly excited states of complex nuclei.
These matrices are taken to have real, complex or real quaternion elements,
and correspond to the quantum system having time reversal symmetry, no
time reversal symmetry, or time reversal symmetry with an odd number
of spin $1/2$ particles respectively.

In the real case all independent elements are chosen from
independent standard normals N$[0,1]$, and in the complex 
case the diagonal elements are chosen independently from N$[0,1]$ 
while the off-diagonal elements are chosen independently
from N$[0,1/\sqrt{2}] + \iup\, \text{N}[0,1/\sqrt{2}]$. 
Real quaternion elements are
themselves $2 \times 2$ blocks of the form
\[
\begin{bmatrix}
\hphantom{-} z & w \\ -\bar{w} & \bar{z} \end{bmatrix}.
\]
In general, the eigenvalues of matrices with real quaternion elements
are doubly degenerate.
On the diagonal, each independent entry is chosen from
N$[0,\sqrt{2}]$, while on the off-diagonal, each independent entry is 
chosen from N$[0,1] + \iup\, \text{N}[0,1]$. 
These ensembles of random matrices are referred to as the Gaussian 
orthogonal, unitary and symplectic ensembles (abbreviated GOE, GUE and GSE)
respectively.

For each of the three Gaussian ensembles the joint probability 
density function (PDF) for the eigenvalues can be computed 
explicitly as \cite{Mehta04}
\begin{equation}\label{Fz}
\frac{1}{(2\pi)^{n/2} F_n(\beta/2)}
\prod_{i=1}^n \eup^{-t_i^2/2} 
\prod_{1 \le i < j \le n}
\abs{t_i - t_j}^{\beta}.
\end{equation}
Here $\beta=1,2,4$ for the GOE, GUE and GSE respectively, while 
$F_n$ is the normalization
\begin{equation}\label{FNz}
F_n(\gamma):=
\frac{1}{(2\pi)^{n/2}}
\int_{-\infty}^{\infty} \cdots \int_{-\infty}^{\infty}
\prod_{i=1}^n \eup^{-t_i^2/2} \prod_{1 \le i < j \le n}
\abs{t_i - t_j}^{2\gamma}\,\dup t_1\cdots\dup t_n,
\end{equation}
referred to as Mehta's integral.
In \cite{MD63} Mehta and Dyson evaluated $F_n(\beta/2)$ for each of the
the three special values of $\beta$. Combining this with
the evaluations for $n=2$ and $n=3$ for general $\beta$ led them
to conjecture that 
\begin{equation}\label{FNz1}
F_n(\gamma) = \prod_{j=1}^n
\frac{\Gamma(1+j\gamma)}{\Gamma(1 + \gamma)}.
\end{equation}

Assuming the validity of \eqref{FNz1} for general $\gamma$, 
certain averages associated with \eqref{Fz} at the special random 
matrix couplings $\beta=1,2,4$ are accessible. This becomes apparent
by writing \eqref{Fz} in the form
\begin{equation}\label{Fza}
\frac{\eup^{-\beta U}}
{(2\pi)^{n/2} F_n(\beta/2)}, 
\qquad\quad U = \frac{1}{2\beta} \sum_{i=1}^n t_i^2  -
\sum_{1 \le i < j \le n} \log \abs{t_i - t_j}. 
\end{equation}
The mean $\langle U \rangle$ 
for a given $\beta$ is now
computed by taking the logarithmic derivative of the
normalization $F_n(\beta/2)$. A further differentiation with respect
to $\beta$ then yields the fluctuation 
$\langle U^2 \rangle - \langle U \rangle^2$.

The form \eqref{Fza} highlights an analogy with the equilibrium 
statistical mechanics of a classical gas of $n$ particles on the line, 
at inverse temperature $\beta$, interacting via a repulsive logarithmic 
Coulomb potential and confined by a harmonic well. 
The quantity $\exp(-\beta U)$ is referred to as the 
Boltzmann factor. This interpretation plays a prominent role in Dyson's series 
of works. 
Indeed, the notation for the averages used 
above stems from the statistical physics literature
(and corresponds to the mean energy and specific heat of the
Coulomb gas) and may be substituted by the mean $\mu(U)$ and variance 
$\sigma^2(U)$ respectively.

It is not hard to see that the Selberg integral can be used 
to evaluate Mehta's integral thus proving the conjecture \eqref{FNz1}.
By the change of variables $t_i \mapsto (1-t_i/L)/2$ in \eqref{SelInt}
\begin{equation}\label{SFlimit}
\lim_{L\to\infty} 2^{L^2} (2L)^{n +n(n-1)\gamma}
S_n(L^2/2,L^2/2,\gamma)=F_n(\gamma).
\end{equation}
Use of Stirling's formula to compute the same limit on the
right-hand side of \eqref{SelInt} then gives \eqref{FNz1}.
However, in 1963 when Mehta and Dyson published their conjecture
the Selberg integral was essentially unknown
and so this method of proof was not available.

The Mehta--Dyson conjecture received more prominence with its
appearance in the first edition of Mehta's book \textit{Random Matrices and 
the Statistical Theory of Energy Levels}, published in 1967 \cite{Mehta67}.
In 1974 Mehta submitted the conjecture to the problems section 
of \textit{SIAM Review} \cite{Mehta74}, 
thus gaining exposure to an even wider
mathematical audience. 
A proof, exactly the one mentioned in the
previous paragraphs, was finally uncovered in the late 1970s by
Selberg's IAS colleague Enrico Bombieri.
The remarkable story behind this proof
is best told in Bombieri's own words \cite{Bombieri07}:
\begin{quote}
Since 1976 I had been studying elementary methods in prime
number theory and in particular a several variable extension of
Chebyshev's well-known method to obtain upper and lower bounds
for the number of primes up to a given bound.  In the course of my
researches I came across the problem of the asymptotic computation
of certain multiple integrals, the simplest being
\[
\int \prod_{i=1}^n z_i^{-r-1}(1-z_i)^p 
\prod_{\substack{i,j=1 \\ i\neq j}}^n (z_i-z_j)^q\, \dup z_1\cdots \dup z_n
\]
where $p,q,r$ are large positive integers and the integral is extended
to the product of the unit circles $\abs{z_i}=1$ or to $[0,1]^n$.

The integral is related to a partition function for the
one-di\-mensio\-nal 
Coulomb gas on the unit circle $\abs{z}=1$ with a fixed point
charge at $z=1$, as it was explained to me by my friend and colleague Tom
Spencer, so I went to Dyson and asked him whether physicists had
encountered such things before; maybe he could save me some efforts.

Dyson told me that for $q=1$ and $2$ an integral of this type, over the real
line with a gaussian measure, had indeed been studied and he referred me
to a book by Mehta. Then I went to see Atle to ask his opinion about
what I was doing in order to study the distribution of primes and
whether he felt it was of any interest and whether he had any opinion
on it.

He immediately recognized my integral as a complex version of the
generalized beta integral he had studied before and he gave me an
off-print of his paper. 
It was not difficult to follow his proof, given for
an integral over $[0,1]^n$, and use a classical method to write a Beta
integral as a complex integral to solve my problem of computing my
integral exactly.  
The multiple integral over
$[0,1]^n$ is of course Selberg's integral, as in that case
arithmetical applications require $r$ to be a large negative
(not positive) integer.
It was also quite easy to get a confluent form of
the Selberg integral and compute exactly the Mehta integrals
for a general value of the parameter and make physicists happy.

Since this was of interest to Dyson, I went back to Dyson and told
him that using the Selberg integral one could compute the integral
of interest to physicists.
\end{quote}

\subsection*{More from the 1960s and 70s --- constant term identities}
The consideration of time reversal symmetry leading to three ensembles of 
Hermitian matrices applies equally well to unitary matrices \cite{Dyson62}.
A conventional time reversal symmetry requires that $U = U^T$, no time reversal
symmetry imposes no constraint, whilst a time reversal symmetry for a 
system with an odd number of spin $1/2$ particles requires $U = U^D$ 
(here $D$ denotes the quaternion dual; see e.g.,~\cite[Ch.~2]{Fo02}). 
Choosing such matrices with a uniform probability then gives what are 
referred to as the circular orthogonal, unitary and symplectic
ensembles (COE, CUE and CSE) respectively. Their joint eigenvalue PDFs are 
given explicitly by
\begin{equation}\label{Cp}
\frac{1}{(2\pi)^n C_n(\beta/2)} \prod_{1 \le i < j \le n}
\abs{\eup^{\iup \theta_i} - \eup^{\iup \theta_j} }^{\beta},
\end{equation}
where $C_n$ is the normalization
\begin{equation}\label{CNb}
C_n(\gamma) = \frac{1}{(2\pi)^n} \int_{-\pi}^{\pi} \cdots \int_{-\pi}^{\pi}
\prod_{1 \le i < j \le n}
\abs{\eup^{\iup \theta_i} - \eup^{\iup \theta_j} }^{2\gamma} \, 
\dup \theta_1 \cdots \dup \theta_n,
\end{equation}
and $\beta = 1,2,4$ for the COE, CUE and CSE respectively.

As for \eqref{FNz}, the random matrix calculations give \eqref{CNb} 
in terms of gamma functions for the three special values of $\beta$. 
Furthermore, the case $n=2$ for general $\beta$ can be related to 
the Euler beta integral \eqref{betaInt}, whilst the case $n=3$ gives
a sum which is a special instancec of an identity of Dixon for a
well-poised $_3F_2$ series \cite{AAR99},
\begin{multline}\label{3F2}
{_3F_2}\biggl(\genfrac{}{}{0pt}{}
{a,b,c}{1+a-b,1+a-c};1\biggr) \\
=\frac{\Gamma(1 + \frac{a}{2}) \Gamma(1 + a - b)  \Gamma(1 + a - c)
\Gamma(1 + \frac{a}{2} - b - c)}{\Gamma(1+a) \Gamma(1 + \frac{a}{2} - b)
\Gamma(1 + \frac{a}{2} - c) \Gamma(1+a-b-c) }.
\end{multline}
Based on all of these results, Dyson, in part I of his series of papers,
made the conjecture \cite{Dyson62}
\begin{equation}\label{CNb1}
C_n(\gamma)=\frac{\Gamma(1+n\gamma)}{\Gamma^n(1+\gamma)}.
\end{equation}
In the same paper, Dyson observed that with $\gamma$ a nonnegative
integer, say $k$, 
\eqref{CNb} can be rewritten as the constant term (CT) in a Laurent 
expansion. This allows \eqref{CNb1} to be rewritten as
\begin{equation}\label{CNb2}
\text{CT} \prod_{1 \le i < j \le n} \Bigl( 1 - \frac{x_i}{x_j} \Bigr)^k
 \Bigl( 1 - \frac{x_j}{x_i} \Bigr)^k = \frac{(kn)!}{(k!)^n}.
\end{equation}
This constant term identity, and thus, by Carlson's theorem,
the conjecture \eqref{CNb1}, was soon proved by J.~Gunson
and K.~Wilson \cite{Wilson62}, and later in a very efficient analysis by 
I.J.~Good \cite{Good70}. Gunson's proof is mentioned in \cite{Dyson62}, but 
the work is unpublished; reference often given to \cite{Gunson62} in this 
context actually refers to the proof of another conjecture of Dyson.
Twenty years after his proof Wilson was to receive the Nobel Prize in
physics for his work on the renormalisation group approach to the
study of critical phenomena; \cite{Wilson62} is his first publication.

In their proof, Wilson and Good both took advantage of the extra degrees 
of freedom offered by Dyson's more general conjecture, 
also contained in \cite{Dyson62},
\begin{equation}\label{CNb3}
\text{CT} \prod_{\substack{i,j=1 \\i \ne j}}^n
\Bigl( 1 - \frac{x_i}{x_j} \Bigr)^{a_i} = 
\frac{(a_1 + \cdots + a_n)!}{a_1! \cdots a_n!}.
\end{equation}
The formulation of this was in turn motivated by the extra 
degrees of freedom permitted by Dixon's identity \eqref{3F2},
to which \eqref{CNb3} reduces in the case $n=3$.

In fact, as observed by R.~Askey \cite{Askey80},
the Selberg integral can be used to prove Dyson's conjecture
\eqref{CNb1} directly without the need for \eqref{CNb3}. 
Askey's observation is based on the easily established general 
identity 
\begin{multline}\label{si}
\int_0^1 \cdots \int_0^1 
(t_1\cdots t_n)^{\zeta-1}
f(t_1,\dots,t_n) \,  \dup t_1 \cdots \dup t_n \\[2mm]
=\Bigl( \frac{1}{2\sin \pi \zeta} \Bigr)^n
\int_{-\pi}^{\pi} \cdots \int_{-\pi}^{\pi}  
\eup^{\iup\zeta(\theta_1+\cdots+\theta_n)}
f(- \eup^{\iup\theta_1} ,\dots, - \eup^{\iup\theta_n}) \, 
\dup \theta_1 \cdots \dup \theta_n,
\end{multline}
valid for $f$ a Laurent polynomial and Re$(\zeta)$ large 
enough so that the left-hand side exists.
Applying \eqref{si} to the Selberg integral with $\beta$ a positive integer
and $\gamma$ a nonnegative integer shows that
\begin{equation}\label{SM}
S_n(\alpha,\beta,\gamma)=(-1)^{n+\binom{n}{2}\gamma}
\Bigl(\frac{\pi}{\sin \pi b}\Bigr)^n M_n(a,b,\gamma),
\end{equation}
where $\alpha:=-b-(n-1)\gamma$, $\beta:=a+b+1$ and
\begin{multline}\label{Morris}
M_n(a,b,\gamma) := \frac{1}{(2\pi)^n}
\int_{-\pi}^{\pi} \cdots \int_{-\pi}^{\pi} 
\prod_{i=1}^n \eup^{\frac{1}{2}\iup\theta_i (a-b)}
\abs{1+\eup^{\iup\theta_i}}^{a+b} \\ \times
 \prod_{1 \le i < j \le n}
\abs{\eup^{\iup \theta_i}-\eup^{\iup \theta_j}}^{2\gamma} \, 
\dup \theta_1 \cdots \dup \theta_n.
\end{multline}
{}From \eqref{SM}, the Selberg integral,
the reflection formula and finally Carlson's theorem,
it follows that
\begin{equation}\label{MorrisEval}
M_n(a,b,\gamma) =
\prod_{j=0}^{n-1} 
\frac{\Gamma (1+a+b+j\gamma) \Gamma(1+(j+1)\gamma)}
{\Gamma (1+a+j\gamma)\Gamma (1+b+j\gamma) \Gamma (1+\gamma)},
\end{equation}
for $a,b,\gamma\in\Complex$ such that 
\[
\Re(a+b+1)>0,~\Re(\gamma)>-\min\{1/n,\Re(a+b+1)/(n-1)\}.
\]
For $a=b=0$ this is Dyson's conjecture \eqref{CNb1}.

\medskip

The change of variables $\eup^{\iup \theta_i}=
(\iup-t_i)/(\iup+t_i)$ maps the unit circle onto the real
line via a stereographic projection. Applying this to the
integral \eqref{Morris} leads to
\begin{multline}\label{SC}
\frac{1}{(2\pi)^n}
\int_{-\infty}^{\infty}\cdots \int_{-\infty}^{\infty}
\prod_{i=1}^n  \frac{1}{(1 + \iup t_i)^{\alpha}
 (1 - \iup t_i)^{\beta}}
\prod_{1 \le i < j \le n} \abs{t_i - t_j}^{2 \gamma} \,
\dup t_1 \cdots \dup t_n \\
=2^{-n(\alpha+\beta-1)+n(n-1)\gamma}
\prod_{j=0}^{n-1}
\frac{\Gamma (\alpha+\beta-1-(n+j-1)\gamma) \Gamma(1+(j+1)\gamma)}
{\Gamma (\alpha-j\gamma)\Gamma (\beta-j\gamma) \Gamma (1+\gamma)}.
\end{multline}
When $n=1$ this is the Cauchy beta integral.

In the letter to Dyson reprinted on the next page
Selberg communicated the multiple Cauchy integral \eqref{SC}. 
Subsequently, in a letter to Askey dated 25 March 1980 \label{pageAskey},
he mentioned both \eqref{Morris} and \eqref{SC},
and pointed out their exact relationship.
The first time \eqref{Morris} appeared in print
was in W.G. Morris' 1982 PhD thesis \cite{Morris82}.
In his thesis Morris provided a proof of \eqref{Morris}
along the lines of Selberg's proof of \eqref{SelInt}, and applied 
it to obtain constant term identities. For these reasons
\eqref{Morris} is now commonly referred to as the Morris integral.

\newpage

\label{pagebrief}
\begin{picture}(0,0)
\put(-65,-650){\includegraphics[width=165mm]{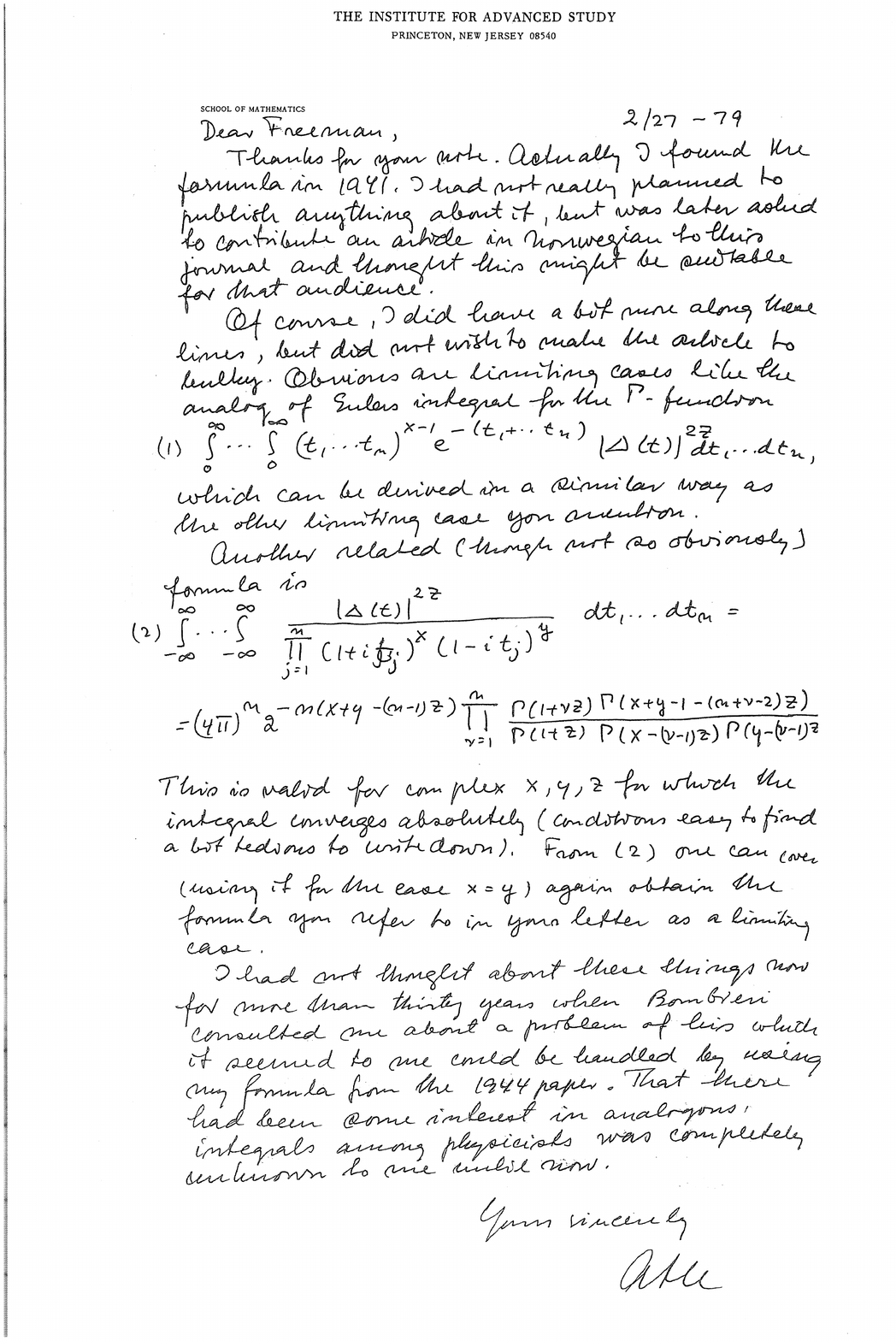}}
\end{picture}\newpage

\newpage

\subsection*{A culmination --- The Macdonald Conjectures}

In 1982 I.G.~Macdonald \cite{Macdonald82} published his now famous
ex-conjectures, generalizing the Mehta integral \eqref{FNz} to all
finite reflection or Coxeter groups $G$, and the Dyson constant
term identity \eqref{CNb2} to all finite root systems.

Let $G$ be a finite group of isometries of $\R^n$,
generated by reflections in $N$ hyperplanes.
Normalise (up to sign) so that the equations for the
hyperplanes take the form
\[
a_{i1} x_1+\cdots+a_{in} x_n=0\quad\text{with}\quad
a_{i1}^2+\cdots+a_{in}^2=2,
\]
where $i$ labels the hyperplanes, and form
the polynomial
\[
P(x)=\prod_{i=1}^N (a_{i1}x_1+\cdots+a_{in} x_n).
\]
Geometrically, $2^{-N/2}P(x)$ gives the product of the
distances of the point $x \in \R$ to the $N$
hyperplanes.

By its action on $\R^n$ the group $G$
acts on polynomials in $x=(x_1,\dots,x_n)$.
The polynomials that are invariant under the action of 
$G$ are referred to as $G$-invariant polynomials.
They form an $\R$-algebra $\R[f_1,\dots,f_n]$ generated by
$n$ algebraically independent polynomials $f_1,\dots,f_n$
of degrees $d_1,\dots,d_n$. Unlike the set of $f_i$'s, the set of
$d_i$'s is uniquely determined by the underlying
reflection group.

A final ingredient required in the Macdonald integral conjectures
is the Gaussian measure $\varphi$ on $\R^n$
\[
\dup\varphi(x):=\frac{\eup^{-\abs{x}^2/2}}{(2\pi)^{n/2}}\,
\dup x_1\cdots\dup x_n,
\]
where $\abs{x}^2:=\sum_{i=1}^n x_i^2$.

With the above notations Macdonald's (ex)-conjecture 
\cite[Conjecture 5.1]{Macdonald82} states that
for each finite reflection group $G$
\begin{equation}\label{ConG}
\Int_{\R^n}
\abs{P(x)}^{2\gamma}\, \dup \varphi(x)=
\prod_{i=1}^n \frac{\Gamma(1+d_i\gamma)}{\Gamma(1+\gamma)}.
\end{equation}
For the three infinite families of crystallographic reflection
groups (or reflection groups of Weyl type) A$_{n-1}$, B$_n$ and D$_n$
the Macdonald conjecture follows as a limit of the Selberg integral.
For type A$_{n-1}$ this corresponds to Bombieri's proof of the
Mehta integral mentioned earlier. 
For types B$_n$ and D$_n$ this is due to A. Regev,
although the actual proof appeared for the first time in
the paper by Macdonald, to whom Regev communicated his results.

Around the same time as Macdonald formulated his conjectures
Regev was studying the large $n$ behaviour of sums of the form
$S_{\ell}^{\beta}(n):=\sum_{\la} (d_{\la})^{\beta}$
where the sum is over partitions $\la$ of at most $\ell$ parts,
and $d_{\la}$ is the dimension of the irreducible $\Symm_n$ 
character $[\la]$. Combining the hook-length formula
for $d_{\la}$ with a careful
asymptotic analysis, Regev showed \cite{Regev81} (see also \cite{CR99})
that the asymptotics of sums like $S_{\ell}^{\beta}(n)$ leads exactly 
to Mehta's integral. Regev remarks \cite{Regev07}
\begin{quote}
{}From reactions to preprints and talks on \cite{Regev81}, first from Richard
Stanley (who in 1978 attended my seminar talk at UCSD) then from
Freeman Dyson, I learned about the Mehta and the Macdonald
conjectures. In a letter, Dyson also mentioned that the Mehta
conjecture had just been verified --- by applying the Selberg
integral. William Beckner then showed me the details of how to
deduce the Mehta --- and some other integrals --- from the Selberg
integral.
I worked on the other classical cases of the Macdonald conjecture
and managed to verify these shortly afterwards, in 1979. 
\end{quote}

The Coxeter group A$_{n-1}$ is the symmetry group of the
$(n-1)$-simplex.
It is a group of order $n!$ generated by the $n(n-1)/2$ hyperplanes
\[
x_i-x_j=0 \quad\text{for}\quad 1\leq i<j\leq n,
\]
and is isomorphic to the symmetric group $\Symm_n$.
All of the ingredients in \eqref{ConG} can thus
easily be determined explicitly.
The polynomial $P(x)$ is given by the Vandermonde product
\begin{equation}\label{VanderM}
P(x)=\prod_{1\leq i<j\leq n} (x_i-x_j) =: \Delta(x)
\end{equation}
and the $G$-invariant polynomials are given by the symmetric
polynomials in $x$. One of the classical bases for the algebra
of symmetric functions is given by the elementary symmetric functions
$\{e_1,\dots,e_n\}$ with
\[
e_r(x)=\sum_{1\leq i_1<i_2<\dots<i_r\leq n} x_{i_1}x_{i_2}\cdots x_{i_r}.
\]
Accordingly the set of degrees $d_i$ is given by
$\{1,2,\dots,n\}$, and \eqref{ConG} reduces to
Mehta's integral \eqref{FNz}. 

The Coxeter groups B$_n$ and D$_n$ are the symmetry groups
of the $n$-cube and $n$-demicube, and for these groups
\eqref{ConG} takes the form
\begin{align*}
\Int_{\R^n}
\prod_{i=1}^n (2\abs{x_i}^2)^{\gamma}
\prod_{1\leq i<j\leq n}\abs{x_i^2-x_j^2}^{2\gamma} \, \dup\varphi(x)
=\prod_{i=1}^n \frac{\Gamma(1+2i\gamma)}{\Gamma(1+\gamma)}, \\
\intertext{and}
\Int_{\R^n}
\prod_{1\leq i<j\leq n}\abs{x_i^2-x_j^2}^{2\gamma} \,
\dup\varphi(x)
=\frac{\Gamma(1+n\gamma)}{\Gamma(1+\gamma)}
\prod_{i=1}^{n-1} \frac{\Gamma(1+2i\gamma)}{\Gamma(1+\gamma)},
\end{align*}
respectively.
Making the changes $t_i=x_i^2/(2L)$, $\alpha=c+1/2$ and $\beta=L+1$
in the Selberg integral and letting $L$ tend to infinity gives 
\[
\Int_{\R^n}
\prod_{i=1}^n (2\abs{x_i}^2)^c
\prod_{1\leq i<j\leq n}\abs{x_i^2-x_j^2}^{2\gamma} \, \dup\varphi(x)
=\prod_{j=0}^{n-1}
\frac{\Gamma(1+2c+2j\gamma)\Gamma(1+(j+1)\gamma)}
{\Gamma(1+c+j\gamma)\Gamma(1+\gamma)},
\]
where on the left use has been made of Legendre's
duplication formula.
The above integral, known as the BC$_n$ Mehta integral,
leads to the B$_n$ and D$_n$ integrals
by setting $c=\gamma$ and $c=0$ respectively.

In his original paper Macdonald established several
other instances of his conjecture, not relying on
the Selberg integral. For $\gamma=1$ Macdonald
presented a uniform proof for all crystallographic
reflection groups. Another case of \eqref{ConG}
--- one that may be verified by purely elementary means ---
is that of the dihedral group I$_2(m)$, the
symmetry group of a regular $m$-gon.

A uniform proof of Macdonald's conjecture for all crystallographic
reflection groups --- A$_{n-1}$, B$_n$, D$_n$, E$_6$, E$_7$, E$_8$, F$_4$
and G$_2$ --- was found in 1989 by E.~Opdam \cite{Opdam89} 
using the Heckman--Opdam theory of hypergeometric shift 
operators \cite{Opdam89,Heckman91}.
Several years later, combined theoretical and computer
efforts by Opdam \cite{Opdam93} and F.~Garvan \cite{Garvan89,Garvan}
also dealt with the remaining non-crystallographic groups $\text{H}_3$ and
$\text{H}_4$, thereby settling Macdonald's conjecture in full.

\medskip

In his paper Macdonald put forward many further conjectures related to
root systems. One of these \cite[Conjecture 2.7]{Macdonald82}
is the generalization of Dyson's constant term identity \eqref{CNb2} 
to arbitrary (finite) root systems.
Let $\Phi$ be a root system (not necessarily reduced) with
corresponding Weyl group $\mathcal W$. For $\alpha\in\Phi$
let $\exp(\alpha)$ be a formal exponential. Denote the
degrees of the fundamental invariants of $\mathcal W$ by $d_1,\dots,d_l$.
The $d_i$ may, for example, be obtained from the simple formula
\[
\prod_{\alpha\in\Phi^{+}}
\frac{1-t^{\height(\alpha)+s(\alpha)}}{1-t^{\height(\alpha)}}
=\prod_{i=1}^l\frac{1-t^{d_i}}{1-t},
\]
where $\Phi^{+}$ is the set of positive roots of the root system,
$\height(\alpha)$ is the height of the root $\alpha$ and $s(\alpha)=1$
if $\alpha/2\not\in\Phi^{+}$ and $s(\alpha)=2$ if $\alpha/2\in\Phi^{+}$.
(The latter can only occur for nonreduced root systems.)
Then Macdonald's constant term conjecture asserts that
\begin{equation}\label{MacCT}
\text{CT}\:\prod_{\alpha\in\Phi}(1-\eup^{\alpha})^k=
\prod_{i=1}^l \binom{d_i k}{k}.
\end{equation}
For the root system A$_{n-1}$,
$\Phi=\{\epsilon_i-\epsilon_j| 1\leq i,j\leq n,~i\neq j\}$
with $\epsilon_i$ the $i$th standard unit vector in $\R^n$.
The degrees in this case are given by $2,3,\dots,n$ so that,
after the identification $\exp(\epsilon_i-\epsilon_j)=x_i/x_j$,
one recovers Dyson's conjecture.

When $k=1$ equation \eqref{MacCT} simply follows from the classical Weyl
denominator formula. Macdonald also proved the $k=2$ case
using algebraic techniques. Once again the Selberg integral  
implies the conjecture for all infinite series: B$_n$, C$_n$, D$_n$ and
BC$_n$. Since the first three are all contained in the latter the
most succinct derivation arises by slightly generalizing the
problem --- Macdonald does this for all root systems in
\cite[Conjecture 2.3]{Macdonald82} --- and considering the 
constant term of
\begin{equation*}
\text{CT}\:\prod_{\alpha\in\Phi_{\text{BC}_n}}
(1-\eup^{\alpha})^{k_{\alpha}}.
\end{equation*}
Here 
\[
\Phi_{\text{BC}_n} = \{\pm \epsilon_i| 1\leq i\leq n\}\cup
\{\pm 2\epsilon_i| 1\leq i\leq n\}\cup
\{\pm \epsilon_i\pm \epsilon_j| 1\leq i<j\leq n\}
\]
is the BC$_n$ root system and $k_{\alpha}=k_1$ if $\alpha$ is of type
$\pm\epsilon_i$, $k_{\alpha}=k_3$ if $\alpha$ is of type $\pm 2\epsilon_i$
and $k_{\alpha}=k_2$ otherwise. The root systems B$_n$, C$_n$ and D$_n$
are then obtained by taking $k_3=0$ or $k_1=0$ or $k_1=k_3=0$ respectively.
By the substitution $\exp(\epsilon_i)\mapsto\exp(2\iup\theta_i)$
it follows that
\begin{multline}\label{BCS}
\text{CT}\:\prod_{\alpha\in\Phi_{\text{BC}_n}}
(1-\eup^{\alpha})^{k_{\alpha}} 
=\frac{2^{n(a+b)+n(n-1)c}}{\pi^n} \\ \times
\int_0^{\pi}\cdots\int_0^{\pi}
\prod_{i=1}^n\sin^a(\theta_i)\cos^b(\theta_i)
\prod_{1\leq i<j\leq n} 
\sin^{c}(\theta_i-\theta_j)
\sin^{c}(\theta_i+\theta_j) \,\dup\theta_1\cdots\dup\theta_n,
\end{multline}
with $a=2k_1+2k_3$, $b=2k_3$ and $c=2k_2$. Introducing new integration
variables $t_i=\sin^2(\theta_i)$ for all $1\leq i\leq n$ 
the integral on the right
transforms into the Selberg integral, so that by \eqref{SelInt} and
the Legendre duplication formula
\begin{align*}
\text{CT}\:\prod_{\alpha\in\Phi_{\text{BC}_n}}&
(1-\eup^{\alpha})^{k_{\alpha}}  \\
&=\frac{4^{n(k_1+2k_3)+n(n-1)k_2}}{\pi^n}\,
S_n\bigl(k_1+k_3+\tfrac{1}{2},k_3+\tfrac{1}{2},k_2\bigr) \\
&=\prod_{i=0}^{n-1} \frac{(k_2+ik_2)!(2k_1+2k_3+2ik_2)!(2k_3+2ik_2)!}
{k_2!(k_1+k_3+ik_2)!(k_3+ik_2)!(k_1+2k_3+(n+i-1)k_2)!}.
\end{align*}
In a not dissimilar manner D. Zeilberger \cite{Zeilberger87} 
showed that the $n=3$ case of the Morris integral 
\eqref{Morris} leads to the Macdonald conjecture
for the exceptional root system G$_2$. This result later found application
in a study linking random matrix theory to number theoretical $L$-functions
\cite{KLR03} (see also the section on the value distribution of
$\log \zeta(1/2+it)$ below).

A unified proof of \eqref{MacCT} for all root systems,
based on hypergeometric shift operators,
is again due to Opdam \cite{Opdam89}. 
Pages \pageref{pageJack1}--\pageref{pageJack2} below
contain an outline of this proof for the root system A$_{n-1}$.

\section{Underpinnings of the Selberg integral}

\subsection*{The Dixon--Anderson integral}

The Euler beta integral \eqref{betaInt} has for its integrand the 
product of power functions $x^{\alpha-1} y^{\beta-1}$ with $y=1-x$.
It is evaluated as a ratio of gamma functions, which in turn are 
integrals over the product of a power function and exponential function.
In the theory of finite fields, the role of power and exponential 
functions are played by multiplicative and additive characters.
These can be used to define the finite field analogue of the gamma
and beta integrals, known as the Gauss and Jacobi sums respectively.
Moreover, these finite field quantities satisfy an analogue of
the beta integral. From Selberg's commentary \cite{Selberg89}, we know that
in the 1940s he investigated finite field analogues of \eqref{SelInt},
and formulated a conjecture which he could prove only for $n=2$.
The existence of such finite field analogues was revealed by Selberg 
in the letter to Askey dated 25 March 1980, referred to on 
page~\pageref{pageAskey}.
Selberg also mentioned this in some colloquium lectures.
A member of the audience on one of these occasions, Ron Evans,
has provided us with the following recollection \cite{Evans07}:
\begin{quote}
Somewhere around 1980, Selberg came to UCSD for a colloquium talk.
Some department members at the UCSD talk were shocked by the subject matter.
They were expecting to hear about his recent work, but instead his entire
talk was on the Selberg integral.
I was fascinated to learn of this integral, and ended up writing several
papers on $q$-analogues and on finite field analogues.  
One of these (published in 1981) formulated $n$-dimensional finite field 
analogues, which I was able to prove for $n=2$.
Selberg had mentioned in his talk that he had finite field analogues for
$n=2$, so I was reluctant to write up my proof.
However, some people who knew Selberg told me that he'd never ever get 
around to publishing his proof, so I took the bold step of asking 
permission to include my proof for $n=2$ with my general conjecture
(with due credit, of course).
He generously wrote back that he didn't mind if I publish a proof of the
``right'' version of the theorem, but that he didn't want to be credited 
with my version, which was too weak!
So I proved the stronger theorem for $n=2$ that he supplied in his letter,
and that led to stronger conjectures for general $n$
(ultimately proved by Anderson \cite{Anderson90}).
\end{quote}
The finite field paper of Evans referred to above is \cite{Evans81}, and
Selberg in his commentary \cite{Selberg89} references this as being his state
of knowledge from the 1940s. In fact, the Anderson paper left open
some of the conjectures from \cite{Evans81} and Evans himself
was able to apply Anderson's approach to provide the remaining proofs
\cite{Evans91}.
For a detailed account of the finite field Selberg integral, we
refer to \cite{AAR99}.

In 1991, motivated by the quest for a proof of the finite field conjecture,
G.W.~Anderson \cite{Anderson91} published a proof of the Selberg integral
based on another multiple integral, namely
\begin{multline}\label{GWA}
\int_X 
\prod_{1 \le i < j \le n} (t_i - t_j) 
\prod_{i=1}^n \prod_{j=1}^{n+1} \abs{t_i - a_j}^{s_j-1}
\, \dup t_1 \cdots \dup t_n \\
= \frac{\prod_{i=1}^{n+1} \Gamma (s_i)}
{\Gamma\bigl(\sum_{i=1}^{n+1} s_i\bigr)}
\prod_{1 \le i < j \le n+1} (a_i - a_j)^{s_i + s_j - 1},
\end{multline}
where $X$ is the domain of integration \label{DApage}
$a_1 >t_1>a_2>t_2>\cdots>t_n>a_{n+1}$,
and $\Re(s_i)>1$ for $1\leq i\leq n+1$.
Anderson's idea was to use \eqref{GWA} to compute in two different 
ways the integral
\begin{multline*}
K(\alpha,\beta,\gamma):= 
\int_{X'} 
\prod_{i=1}^{n+1} x_i^{\alpha-1} (1 - x_i)^{\beta-1} 
\prod_{i=1}^n \prod_{j=1}^{n+1} \abs{y_i - x_j}^{\gamma-1} \\
\times \prod_{1\leq i<j\leq n} \abs{y_i - y_j} 
\prod_{1\leq i<j\leq n+1} \abs{x_i - x_j}
\, \dup x_1 \cdots \dup x_{n+1} \dup y_1 \cdots \dup y_n,
\end{multline*}
where $X'$ denotes the domain of integration
\[
1 > x_1 > y_1 > x_2 > y_2 > \cdots > y_n > x_{n+1} > 0.
\]
First integrating over the $y$-variables gives
\[
K(\alpha,\beta,\gamma) = 
\frac{\Gamma^{n+1}(\gamma)}{(n+1)!\,\Gamma((n+1)\gamma)}\:
S_{n+1}(\alpha,\beta,\gamma),
\]
while first integrating over the $x$-variables gives
\[
K(\alpha,\beta,\gamma) =
\frac{\Gamma (\alpha) \Gamma (\beta) \Gamma^n(\gamma)}
{n!\,\Gamma (\alpha + \beta +n\gamma)} \:
S_n(\alpha+\gamma,\beta+\gamma,\gamma).
\]
Equating the two forms reveals a first order recurrence for 
the Selberg integral in $n$. Together with the initial condition
$S_0(\alpha,\beta,\gamma)=1$ this reclaims \eqref{SelInt}.

A large portion of Anderson's paper is devoted to a derivation of
\eqref{GWA}. This same multiple integral, written in the form
\begin{equation}\label{V}
\det_{1\leq i,j\leq n}
\biggl( \int_{a_{i+1}}^{a_i} t^{j-1} 
\prod_{l=1}^{n+1} \abs{t - a_l}^{s_l - 1} \, \dup t \biggr)
\end{equation}
was evaluated at around the same time by A.~Varchenko
\cite{Varchenko89,Varchenko90} in his work on hyperplane arrangements.
That \eqref{V} is equal to the integral in \eqref{GWA} is a simple 
consequence of the Vandermonde determinant, a fact made explicit 
in \cite{RZ02}.

Remarkably, in a 1998 paper by M.C.~Berg\`ere \cite{Bergere98} proving a 
conjecture from the theory of Calogero--Sutherland models 
(see page~\pageref{CSpage}) reference is made to \eqref{GWA}, 
citing neither Anderson nor Varchenko, but a paper of 
A.L.~Dixon \cite{Dixon05} published in 1905! 
Indeed, consulting \cite{Dixon05} reveals both \eqref{GWA} ---
obtained via essentially the same analysis as that used in \cite{Anderson91} 
--- and its equivalent determinant form \eqref{V}.

\medskip

A study by Forrester and E.M.~Rains \cite{FR02} provides additional links 
between the Selberg and the Dixon--Anderson integrals.
These apply at the level of the corresponding normalized 
integrands, referred to as the Selberg and Dixon--Anderson densities.
The former will be denoted by $S_n(\alpha,\beta,\gamma;t)$ for
$t=(t_1,\dots,t_n)$.

The first point of note is that the computations of 
Dixon and Anderson can be interpreted as giving the density 
of zeros of the random rational function
\[
R_{n+1}(x):=\sum_{i=1}^{n+1} \frac{w_i}{a_i-x},
\]
where the $w_i$ are distributed according to the Dirichlet distribution
--- to be denoted $D_{n+1}[s_1,\dots,s_{n+1}]$ ---
\[
\frac{\Gamma(s_1 + \cdots + s_{n+1})}{\Gamma(s_1) \cdots \Gamma(s_{n+1})}
\prod_{i=1}^{n+1} w_i^{s_i - 1},
\]
with $w_1,\dots, w_{n+1} > 0$ such that $w_1+\cdots+w_{n+1}=1$.
Motivated by this interpretation, a family of random polynomials 
$A_j(x)$, $1\leq j\leq n$ were defined in \cite{FR02} such that the zeros of 
$A_j(x)$ have PDF $S_j(\alpha_j,\beta_j,\gamma;t)$ with
$\alpha_j:=(n-j)\gamma+\alpha$, $\beta_j:=(n-j)\gamma+\beta$.
Setting $A_{-1}(x) := 0$, $A_0(x) := 1$, and specifying that
$(w_0^{(j)},w_1^{(j)},w_2^{(j)})$ be distributed according to 
the Dirichlet distribution $D_3[\beta_j,(j-1)\gamma,\alpha_j]$,
the polynomials $A_j(x)$ are determined by the random 
three-term recurrence
\begin{equation}\label{A3}
A_j(x) = w_2^{(j)} (x-1) A_{j-1}(x) + 
w_0^{(j)} x A_{j-1}(x) + w_1^{(j)} x (x - 1) A_{j-2}(x).
\end{equation}

Let $\alpha\mapsto\gamma\alpha/2+1$,
$\beta\mapsto\gamma\beta/2+1$, and let the integrand of the
Selberg integral be written in the form $\exp(-2 \gamma U)$ so that
\[
U = - \frac{\alpha}{2} \sum_{i=1}^n \log t_i -
\frac{\beta}{2} \sum_{i=1}^n \log\abs{1 - t_i} -
\sum_{1 \le i < j \le n} \log\abs{t_i-t_j}.
\]
Then in the limit $\gamma \to \infty$ the Selberg density
crystallizes to the minimum of $U$ subject to the constraint
that $0 < t_i < 1$ for each $t_i$. According to a classical result of
Stieltjes (see e.g., \cite{AAR99}) this minimum is unique and occurs at
the zeros of the Jacobi polynomial $P_n^{(\alpha-1,\beta-1)}(x)$.
Indeed in the same limit the three term recurrence \eqref{A3} is no 
longer random, and 
has solution $A_j(x) = \tilde{P}_j^{(n-j+\alpha-1,n-j+\beta-1)}(x)$ with
$\tilde{P}^{(a,b)}(x)$ the Jacobi polynomial normalized
to be monic \cite{FR02}.

The change of variables and limiting procedure giving rise
to \eqref{SFlimit} reduces the Selberg density to the PDF \eqref{Fz}.
The Dixon--Anderson density permits a similar limit, and applied 
with $n \mapsto n+1$, $a_0=0$ and $a_1=1$ this results in the PDF on
$\{t_j\}_{j=1}^{n+1}$ given by
\begin{multline}\label{3.st1}
\frac{1}{\sqrt{2 \pi}\, \Gamma(s_1) \cdots \Gamma(s_n)}\,
\frac{\prod_{1 \le i < j \le n+1} (t_i -t_j)}
{\prod_{1 \le i < j \le n} (a_i - a_j)^{s_i+s_j-1}} \\
\times \prod_{i=1}^{n+1} \prod_{j=1}^n \abs{t_i - a_j}^{s_j-1}
\exp \biggl( -\frac{1}{2} \sum_{j=1}^{n+1} t_j^2 +\frac{1}{2}
\sum_{j=1}^n a_j^2 \biggr)
\end{multline}
supported on
\[
t_1 > a_1 > t_2 > \cdots > a_n > t_{n+1}.
\]
The corresponding limit of $R_{n+1}(x)$ gives the random 
rational function
\begin{equation}\label{Rtilde}
\tilde{R}_{n+1}(x) : =
x - \mu_0 + \sum_{i=1}^n \frac{\mu_i}{a_i-x},
\end{equation}
where $\mu_0$ has distribution N$[0,1]$
and $\mu_i$ the Gamma distribution $\Gamma[s_i,1]$.
Indeed, the fact that the zeros of \eqref{Rtilde} have PDF \eqref{3.st1} 
can readily be checked directly by adopting
the strategy of Dixon and Anderson, see \cite{Evans94,Fo06}.
Finally, the limiting form of the three term recurrence \eqref{A3} is
\begin{equation}\label{Cs}
C_j(x) = (x - r) C_{j-1}(x) - s^{(j-1)} C_{j-2}(x)
\end{equation}
with $C_{-1}(x) = 0$, $C_0(x) = 1$, $r$ having distribution N$[0,1]$
and $s^{(j)}$ distribution $\Gamma[j\gamma,1]$. The random polynomial
$C_j(x)$ has as the PDF for its zeros the density \eqref{Fz} with $n=j$
and $\beta=2\gamma$.

It should be remarked that since \eqref{3.st1} integrated over
$t_1,\dots,t_{n+1}$ gives unity, a limiting form of the Dixon--Anderson
integral follows. Evans \cite{Evans94} used this, together with the
strategy of Anderson, to give the first proof of the Mehta integral
evaluation \eqref{FNz1} which is independent of the Selberg integral.

The random polynomial $C_{n+1}(x)$ can be interpreted as the 
characteristic polynomial for a family of random matrices 
defined inductively by \cite{FR02}
\begin{equation}\label{MM}
M_{n+1} =
\begin{bmatrix}
\la_1^{(n)} &\!\! &\!\! & b_1 \\
&\!\!\ddots &\!\! & \vdots \\
&\!\! & \!\! \la_n^{(n)} & b_n \\
b_1 &\!\! \dots &\!\! b_n & c 
\end{bmatrix}.
\end{equation}
Here the $\la_i^{(n)}$ are the eigenvalues of $M_n$,
$c$ has distribution N[0,1] and $b_j^2$ has distribution 
$\Gamma[\beta/2,1]$. Indeed it is straightforward to show that
the eigenvalues of \eqref{MM} are given by the zeros of \eqref{Rtilde},
with $\mu_0=c$, $\mu_i = b_i$ and $a_i=\la_i^{(n)}$. 
In the case $\beta = 1$, the invariance of members of the GOE with respect 
to conjugation by orthogonal matrices shows that \eqref{MM} is similar to 
GOE matrices, and an analogous understanding of the
relationship between \eqref{MM} in the case $\beta=2$ and 
GUE matrices can be given.
Moreover, it is generally true that a three-term recurrence
\[
p_j(x) = (x - a_j) p_{j-1}(x) - b_{j-1} p_{j-2}(x)
\]
with $p_{-1}(x) := 0$, $p_0(x) = 1$ is satisfied by the 
characteristic polynomial for the tridiagonal matrix
\[
\begin{bmatrix}
a_n & b_{n-1} & & & \\
b_{n-1} & a_{n-1} & b_{n-2} & & \\
& \ddots & \ddots & \ddots & \\
& &b_2 & a_2 & b_1 \\
& & & b_1 & a_1 
\end{bmatrix}.
\]
Hence $C_n(x)$ is also the characteristic polynomial for the 
above random tridiagonal matrix with each $a_j$ having 
distribution N[0,1] and with $b_j^2$ distributed as in \eqref{MM}.
This is a result due to I.~Dumitriu and A.~Edelman \cite{DE02}, obtained
without the use of \eqref{3.st1}. 
In this regard it should be mentioned that R.~Killip and I.~Nenciu \cite{KN04}, 
in a study which does not make use of the Dixon--Anderson integral
\eqref{GWA}, give the explicit 
construction of a family of
random orthogonal matrices with eigenvalue PDF equal to the 
BC$_n$ Selberg density, which itself is proportional to the integrand in
\eqref{BCS}. The methods of \cite{DE02} and  \cite{KN04}, which at
a technical level proceed via a change of variables from a general tridiagonal
matrix and unitary Hessenberg matrix to their eigenvalue/eigenvector 
decomposition, yield too the evaluations of the Mehta and Selberg 
integrals respectively.

\subsection*{Dotsenko--Fateev integrals}

In the course of studies in conformal field theory,
V.S.~Dotsenko and V.A.~Fateev \cite{DF85} were lead to 
consider the multiple integral
\begin{multline*}
\text{PV}
\Int_{[0,1]^p}\Int_{[1,\infty)^{n-p}} 
\Int_{[0,1]^r}\Int_{[1,\infty)^{m-r}} 
\prod_{i=1}^n t_i^{\alpha}(1-t_i)^{\beta}
\prod_{i=1}^m \tau_i^{\alpha'} (\tau_i-1)^{\beta'} 
\prod_{i=1}^n \prod_{j=1}^m (\tau_j-t_i)^{-2} \\
\times
\prod_{1\leq i<j\leq n} \abs{t_i-t_j}^{2\gamma}
\prod_{1\leq i<j\leq m} \abs{\tau_i-\tau_j}^{2\gamma'}\,
\dup t_1 \cdots \dup t_n \, \dup \tau_1 \cdots \dup\tau_m,
\end{multline*}
where PV denotes the principal value,
$\alpha/\alpha'=\beta/\beta'=-\gamma$, $\gamma\gamma'=1$ 
and $0\leq p\leq n$, $0\leq r\leq m$.
Note that the case $p=n$ and $m=0$ is, 
up to a shift by $1$ in $\alpha$ and $\beta$,
precisely the Selberg integral. 
Dotsenko and Fateev evaluated the above integral as a product of 
gamma and sine functions reclaiming the Selberg integral as a special case.

The method employed by Dotsenko and Fateev for evaluating their
integral suggests an approach \cite{Fo02} to the Selberg integral 
by studying the simpler $m=0$ case
\begin{equation}\label{DF}
S_{n,p}(\alpha,\beta,\gamma) := 
\Int_{[0,1]^p} \Int_{[1,\infty)^{n-p}}
\prod_{i=1}^n t_i^{\alpha - 1} \abs{1 - t_i}^{\beta - 1} 
\prod_{1 \le i < j \le n} \abs{t_i - t_j}^{2 \gamma}\,
\dup t_1 \cdots \dup t_n
\end{equation}
for $0\leq p\leq n$. Note that 
$S_{n,n}(\alpha,\beta,\gamma)=S_n(\alpha,\beta,\gamma)$, which is the
Selberg integral. Also note that the change of variables 
$t_i \mapsto 1/t_i$ for all $1\leq i\leq n$ implies the transformation
\begin{equation}\label{trafo}
S_{n,p}(\alpha,\beta,\gamma)=
S_{n,n-p}(1-\alpha-\beta-2(n-1)\gamma,\beta,\gamma).
\end{equation}

Singling out the integration variable $t_p$, viewing the integrand as an
analytic function and replacing the 
interval $[0,1]$ by a contour along a positively oriented, indented 
semi-circle of infinite radius (with indentations at the branch points
$t_p=0,1,t_1,\dots, t_{p-1},t_{p+1},\dots,t_n$) yields the recurrence
\begin{multline}\label{DF1}
S_{n,p}(\alpha,\beta,\gamma) =\frac{p}{n-p+1}\,
\frac{\sin \pi (n-p+1)\gamma  \, \sin \pi (\alpha + \beta
+ (n + p - 2)\gamma )}{\sin \pi p \gamma  \,
\sin \pi (\alpha + (p-1) \gamma )} \\[2mm]
\times S_{n,p-1}(\alpha,\beta,\gamma).
\end{multline}
Solving $S_{n,n}$ in terms of $S_{n,0}$ and then using the 
transformation \eqref{trafo} to eliminate $S_{n,0}$ in favour
of $S_{n,n}$ gives the following
functional equation for the Selberg integral
\[
S_n(\alpha,\beta,\gamma) =
S_n(1-\alpha-\beta-2(n-1)\gamma,\beta,\gamma)
\prod_{j=0}^{n-1} \frac{\sin\pi(\alpha +\beta+(n+j-1)\gamma)}
{\sin\pi(\alpha+j\gamma)}.
\]
The significance of this result is that it permits $S_n$, viewed
as a function of $\alpha$, to be analytically continued to 
$\alpha\in\Complex$ with the exclusion of the zeros of
$\sin\pi(\alpha+j\gamma)$. 
Indeed, for sufficiently small but positive 
values of $\Re(\beta)$ and $(n-1)\Re(\gamma)$ the Selberg integral 
requires $\Re(\alpha)>0$ but the transformed integral permits 
$\Re(\alpha)<1-\Re(\beta)-2(n-1)\Re(\gamma)$ which is greater than $0$.
As a consequence one can verify that $S_n(\alpha,\beta,\gamma)$ when divided
by the right-hand side of \eqref{SelInt} is a bounded analytic function 
in the complex $\alpha$-plane, and is thus independent of $\alpha$. 
Symmetry then gives that this ratio is independent of $\beta$ as well.
That the dependence on $\gamma$ and $n$ is correct is then
verified by using
\[
\lim_{\alpha \to 0^{+}} \alpha S_n(\alpha,\beta,\gamma) = 
n S_{n-1}(2 \gamma - 1,\beta,\gamma),
\]
a fact already noted and used for the same purpose in Selberg's 
original proof \cite{Selberg44}.

\medskip 

In their paper Dotsenko and Fateev considered a further 
generalization of \eqref{SelInt}, referred to as the complex Selberg integral.
This integral, which was also studied independently by 
K.~Aomoto \cite{Aomoto87b}, can be written as an $n$-dimensional 
real integral with integration variables given by $2$-dimensional vectors 
\[
A_n(\alpha,\beta,\gamma):=\int_{\R^2} \cdots \int_{\R^2} 
\prod_{i=1}^n \abs{\vec{r}_i}^{2(\alpha-1)}
\abs{\vec{u}-\vec{r}_i}^{2(\beta-1)} 
\prod_{1\le i<j\le n}
\abs{\vec{r}_i-\vec{r}_j}^{4\gamma}\,
\dup\vec{r}_1\cdots\dup\vec{r}_n,
\]
where $\vec{u}$ is an arbitrary unity vector. 
Dotsenko and Fateev as well as Aomoto showed that
up to a product of trigonometric functions the complex Selberg integral
factors as a product of two Selberg integrals
\[
A_n(\alpha,\beta,\gamma)=S_n^2(\alpha,\beta,\gamma) \,
\frac{1}{n!} \prod_{j=0}^{n-1} \frac{\sin \pi(\alpha+j\gamma)
\sin\pi(\beta+j\gamma) \sin\pi (j+1)\gamma}
{\sin \pi (\alpha+\beta+(n+j-1)\gamma) \sin\pi\gamma},
\]
provided \eqref{SelIntpv} is supplemented by
\[
\text{Re}(\alpha+\beta+(n-1)\gamma) < 1 \quad
\text{and}\quad
\text{Re}(\alpha+\beta+2(n-1)\gamma) < 1.
\]

K.~Mimachi and M.~Yoshida \cite{MY04} (see also \cite{MY03})
apply the theory of intersection numbers of 
twisted cycles to the conformal field theory study of Dotsenko and Fateev
to give the evaluation of the product 
$S_n(\alpha,\beta,\gamma)S_n(-\alpha,-\beta,-\gamma)$,
appropriately analytically continued. This is achieved without
requiring the actual evaluation of the Selberg integral itself.

\subsection*{Jack polynomial theory}  

It has been known since the early 1970s \cite{Calogero69,Sutherland71,
Sutherland72} that \eqref{Cp} with $\beta=2\gamma$ ---
to be denoted $\exp(-2\gamma W)$ in analogy with
\eqref{Fza} --- is the absolute value squared of the ground-state 
wave function for the Schr\"odinger operator \label{CSpage}
\[
H = - \sum_{i=1}^n \frac{\partial^2}{\partial \theta_i^2} + 
\frac{1}{2}\gamma(\gamma-1)\sum_{1 \le i < j \le n}
\frac{1}{\sin^2\bigl(\tfrac{1}{2}(\theta_i-\theta_j)\bigr)}.
\]
This operator, known as the Calogero--Sutherland Hamiltonian,
describes a system of $n$ identical quantum particles on the unit circle,
with $\theta_i\in[0,2\pi)$ for $1\leq i\leq n$ the (angular) 
positions of the particles. The interaction between the particles
is described by a $1/r^2$ two-body potential,
$2\abs{\sin((\theta_i - \theta_j)/2)}$
being the cord-length between particles located at
$\theta_i$ and $\theta_j$.

B.~Sutherland \cite{Sutherland71} showed that the eigenvalue $E_0$
corresponding to the ground-state wave function
is given by $E_0=n(n^2-1)\gamma^2/12$. Subsequently he showed 
\cite{Sutherland72} that the conjugated operator
\begin{equation}\label{CO}
\eup^{\gamma W} (H - E_0) \eup^{-\gamma W} =
\sum_{i=1}^n \Bigl(x_i \frac{\partial}{\partial x_i} \Bigr)^2
+ 2\gamma \sum_{\substack{i,j=1 \\ i\neq j}}^n
\frac{x_i + x_j}{x_i - x_j} \, \frac{\partial}{\partial x_i},
\end{equation}
where $x_j := \exp(\iup \theta_j)$,
admits a complete set of symmetric polynomial eigenfunctions
$P_{\la}^{(1/\gamma)}(x)$. These polynomials,
now referred to as Jack polynomials, depend on $x=(x_1,\dots,x_n)$
and are indexed by partitions $\la$ of at most $n$ parts;
$\la=(\la_1,\dots,\la_n)$ with
$\la_1\geq \la_2\geq \dots\geq \la_n\geq 0$. 
With $m_{\la}$ denoting the monomial symmetric 
function indexed by $\la$ and $<$  
the dominance ordering on partitions, the 
Jack polynomials have the structure
\begin{equation}\label{stu}
P_{\la}^{(1/\gamma)}(x) = m_{\la}(x) + 
\sum_{\mu < \la} a_{\la\mu}\, m_{\mu}(x)
\end{equation}
for some coefficients $a_{\la\mu}=a_{\la\mu}(\gamma)$.

One fundamental property of the Jack polynomials is that they 
are orthogonal with respect to the inner product \label{pageJack1}
\begin{equation}\label{InProd}
\langle f,g \rangle_{\gamma}:=\frac{1}{(2\pi)^n}
\int_{-\pi}^{\pi} \cdots \int_{-\pi}^{\pi}
f(\eup^{\iup \theta}) g(\eup^{-\iup \theta})
\prod_{1\le i<j\le n}
\abs{\eup^{\iup\theta_i}-\eup^{\iup\theta_j}}^{2\gamma} \,
\dup \theta_1 \cdots \dup \theta_n,
\end{equation}
where $f(\eup^{\iup \theta})=
f(\eup^{\iup \theta_1},\dots,\eup^{\iup \theta_n})$.
To state the orthogonality as well as the quadratic norm evaluation let
\[
\Poch{b}{\la}{\gamma}=\prod_{i\geq 1} (b+(1-i)\gamma)_{\la_i}
\]
with $(b)_n=b(b+1)\cdots(b+n-1)$ a Pochhammer symbol. Also let
$c_{\la}(\gamma)$ and $c_{\la}'(\gamma)$ be given by
\begin{subequations}\label{ccp}
\begin{align}
c_{\la}(\gamma)&=\prod_{s\in\la}(a(s)+l(s)\gamma+\gamma), \\
c'_{\la}(\gamma)&=\prod_{s\in\la}(a(s)+l(s)\gamma+1),
\end{align}
\end{subequations}
where $a(s)$ and $l(s)$ are the arm-length and leg-length of the 
square $s$ in the diagram of the partition $\la$, and $\abs{\la}$ 
is the total number of boxes in the diagram of $\la$ \cite{Macdonald95}.
Then
\begin{equation}\label{OP}
\bigl\langle P_{\la}^{(1/\gamma)},P_{\mu}^{(1/\gamma)}\bigr\rangle_{\gamma}
=\delta_{\la\mu} \:
\frac{c'_{\la}(\gamma)}
{\Poch{1+(n-1)\gamma}{\la}{\gamma}}\:
\frac{\Gamma(1+n\gamma)}{\Gamma^n(1+\gamma)}\:
P_{\la}^{(1/\gamma)}(1^n),
\end{equation}
where $\delta_{\la\mu}$ is the Kronecker delta function and
$(1^n)$ is shorthand for $(1,1,\dots,1)$.
The orthogonality relation is consistent with, 
but not an immediate consequence of the operator 
\eqref{CO} being self-adjoint with respect to \eqref{InProd}. 
The complication is that not all eigenvalues of \eqref{CO} are distinct.
This degeneracy can be removed by introducing the mutually commuting
Cherednik operators
$\xi_i$ for $1\leq i\leq n$ \cite{Cherednik91,Dunkl89} 
\[
\xi_i=1-i+\frac{x_i}{\gamma} \frac{\partial}{\partial x_i} +
\sum_{j=1}^{i-1}\frac{x_i}{x_i - x_j}\,(1 - s_{ij}) +
\sum_{j=i+1}^n\frac{x_j}{x_i - x_j}\,(1 - s_{ij}),
\]
where $s_{ij}$ acts by permutation 
$x_i$ and $x_j$ and $1$ represents the identity operator.
Any symmetric combination of the $\xi_i$, and in particular
$\prod_{i=1}^n (1-u\xi_i)$,
has the Jack polynomials as simultaneous eigenfunctions.

The Cherednik operators can be used to construct the Jack polynomial 
shift operator --- a special case of the shift operators studied by 
Heckman and Opdam, and used by the latter to prove the Macdonald
integral and constant term conjectures. Properties of the
Jack shift operator not only imply \eqref{CNb1} or, equivalently,
\eqref{CNb2}, but also the more general quadratic norm evaluation 
of the Jack polynomials corresponding to \eqref{OP} 
with $\la=\mu$ \cite{Kakei98}.
(For $\la=0$ this yields \eqref{CNb1}.)
With $\Delta(x)$ the Vandermonde product \eqref{VanderM} and $Y_{\pm}:=
\gamma^{n(n-1)/2} \prod_{1\leq i<j\leq n}(\xi_i-\xi_j\mp 1)$,
the Jack shift operators are defined by 
$G_{+}:=\Delta^{-1} Y_{+}$, $G_{-} = Y_{-}\Delta$.
They have an adjoint type property with respect to the inner product 
\eqref{InProd},
\begin{equation}\label{fr.1}
\langle G_{+}f,g\rangle_{\gamma+1}=\langle f,G_{-} g\rangle_{\gamma}.
\end{equation}
Also, with
\begin{equation}\label{apm}
a_{\la}^{\pm}(\gamma)=\prod_{1\leq i<j\leq n}
(\la_i-\la_j\pm 1+(j-i\mp 1)\gamma)
\end{equation}
and $\delta$ the staircase partition $(n-1,n-2,\dots,1,0)$, 
the shift operators act on the Jack polynomials as
\begin{subequations}\label{fr.2}
\begin{align}
G_{+} P_{\la+\delta}^{(1/\gamma)}&=a_{\la}^{+}(\gamma+1) 
P_{\la}^{(1/(\gamma+1))}, \\[2mm]
G_{-} P_{\la}^{(1/(\gamma+1))}&=a_{\la}^{-}(\gamma+1) 
P_{\la+\delta}^{(1/\gamma )}.
\end{align}
\end{subequations}
It follows from \eqref{fr.1} and \eqref{fr.2} that 
\[
\bigl\langle P_{\la}^{(1/(\gamma+1))},P_{\la}^{(1/(\gamma+1))}
\bigr\rangle_{\gamma+1}=
\frac{a_{\la}^{-}(\gamma+1)}{a_{\la}^{+}(\gamma+1)}\:
\bigl\langle P_{\la+\delta}^{(1/\gamma)},P_{\la+\delta}^{(1/\gamma)}
\bigr\rangle_{\gamma} 
\]
and thus
\[
\bigl\langle P_{\la}^{(1/(\gamma+k))},P_{\la}^{(1/(\gamma+k))}
\bigr\rangle_{\gamma+k}=
\bigl\langle P_{\la+k\delta}^{(1/\gamma)},P_{\la+k\delta}^{(1/\gamma)}
\bigr\rangle_{\gamma}\:
\prod_{j=1}^{k-1} 
\frac{a_{\la+j\delta}^{-}(\gamma+k-j)}{a_{\la+j\delta}^{+}(\gamma+k-j)}.
\]
Taking $\gamma=0$, using that $P_{\la}^{(\infty)}=m_{\la}$ 
(the monomial symmetric function) and  
\[
\bigl\langle m_{\mu},m_{\mu} \bigr\rangle_0=
\text{CT}\Bigl( m_{\mu}(x)m_{\mu}(x^{-1})\Bigr)=
m_{\mu}(1^n)
\]
which is $n!$ for $\mu=\la+k\delta$,
it follows that for nonnegative integer $k$
\label{pageJack2}
\begin{equation}\label{fr.3}
\bigl\langle P_{\la}^{(1/k)},P_{\la}^{(1/k)}\bigr\rangle_k
=n!\prod_{j=0}^{k-1} \frac{a_{\lambda+jk}^{-}(k-j)}
{a_{\lambda+jk}^{+}(k-j)}.
\end{equation}
Using the evaluation formula \cite{Stanley89} 
\begin{equation}\label{ef}
P_{\la}^{(1/\gamma)}(1^n)=\frac{[n\gamma]_{\la}^{(\gamma)}}
{c_{\la}(\gamma)}
\end{equation}
and the definitions \eqref{ccp} and \eqref{apm} it is now
a straightforward exercise to verify that for 
$\gamma=k$ \eqref{OP} coincides with \eqref{fr.3}.
Analytic continuation off the integers is then required to
establish
\eqref{OP} for all $\Re(\gamma)>-1/n$.

\medskip

A further fundamental property of the Jack polynomials is
R.P.~Stanley's Cauchy identity \cite{Stanley89}
\begin{equation}\label{CaP}
\sum_{\la} \frac{c_{\la}(\gamma)}{c'_{\la}(\gamma)} \,
P_{\la}^{(1/\gamma)}(x) P_{\la}^{(1/\gamma)}(y)=
\prod_{i=1}^n\prod_{j=1}^m (1-x_i y_j)^{-\gamma},
\end{equation}
where $x=(x_1,\dots,x_n)$, $y=(y_1,\dots,y_m)$.
Setting $y=(1^m)$, using the evaluation formula 
\eqref{ef} with $n\mapsto m$
and a standard analytic argument to replace $m$ by $a\alpha$,
leads to Z.~Yan's binomial theorem for Jack polynomials \cite{Yan92}
\begin{equation}\label{BTP}
\sum_{\la}
\frac{\Poch{a}{\la}{\alpha}}{c'_{\la}(\alpha)}\,
P_{\la}^{(\alpha)}(x)
=\prod_{i=1}^n \frac{1}{(1-x_i)^a}.
\end{equation}
This, together with the orthogonality \eqref{OP}, the property
\[
P_{(\la_1+a,\dots,\la_n+a)}^{(\alpha)}(x)=
(x_1\cdots x_n)^a P_{\la}^{(\alpha)}(x), \qquad a=0,1,2,\dots
\]
and the gamma reflection formula, implies a
generalization of the Morris integral \eqref{Morris} \cite{BF98}
\begin{multline*}
\frac{1}{(2\pi)^n}
\int_{-\pi}^{\pi}\cdots \int_{-\pi}^{\pi}
P_{\la}^{(1/\gamma)}(-\eup^{\iup\theta})
\prod_{i=1}^n \eup^{\frac{1}{2}\iup \theta_i (a-b)} 
\abs{1+\eup^{\iup\theta_i}}^{a+b} \\ \times
\prod_{1\le i<j\le n}
\abs{\eup^{\iup\theta_i}-\eup^{\iup\theta_j}}^{2\gamma}\,
\dup \theta_1 \cdots \dup \theta_n  \\ 
=\frac{\Poch{-b}{\la}{\gamma}}{\Poch{1+a+(n-1)\gamma}{\la}{\gamma}}\:
P_{\la}^{(1/\gamma)}(1^n) \, M_n(a,b,\gamma).
\end{multline*}
Applying \eqref{si} finally results in a generalization of the 
Selberg integral
\begin{multline}\label{KadellInt}
\int^1_0 \cdots \int^1_0 
P_{\la}^{(1/\gamma)}(t)
\prod_{i=1}^n t_i^{\alpha - 1} (1-t_i)^{\beta - 1}
\prod_{1\leq i<j\leq n}\abs{t_i-t_j}^{2\gamma}\,
\dup t_1 \cdots  \dup t_n \\
=\frac{\Poch{\alpha +(n-1)\gamma}{\la}{\gamma}}
{\Poch{\alpha+\beta+2(n-1)\gamma}{\la}{\gamma}}\:
P_{\la}^{(1/\gamma)}(1^n) \,
S_n(\alpha,\beta,\gamma).
\end{multline}
This evaluation is usually referred to as Kadell's integral \cite{Kadell97}
after its first prover, but as a conjecture is due to Macdonald 
\cite[Conjecture (C5)]{Macdonald87}.
When $\la=(1^r)$, in which case the Jack polynomial
is nothing but the $r$th elementary symmetric function, the above
is known as Aomoto's integral, who used it to give
what is arguably the first elementary proof of the Selberg 
integral \cite{Aomoto87}.
A proof of Kadell's integral along the lines of Anderson's proof 
of the Selberg integral as described on page \pageref{DApage}
has recently been obtained in \cite{W07} through use of the
Okounkov--Olshanski integral formula for Jack polynomials 
\cite{Okounkov98}
\begin{multline*}
P_{\la}^{(1/\gamma)}(x)=
\prod_{i=1}^{n-1}\frac{\Gamma(\la_i+(n-i+1)\gamma)}
{\Gamma(\la_i+(n-i)\gamma)\Gamma(\gamma)}
\prod_{1\leq i<j\leq n}(x_j-x_i)^{1-2\gamma} \\
\times \int_Y
P_{\la}^{(1/\gamma)}(y)\prod_{1\leq i<j\leq n-1}(y_j-y_i)
\prod_{i=1}^{n-1}\prod_{j=1}^n \abs{y_i-x_j}^{\gamma-1} \dup y_1\cdots
\dup y_{n-1}.
\end{multline*}
where $Y$ denotes the domain
$x_1<y_1<x_2<\cdots<x_{n-1}<y_{n-1}<x_n$
and $\la$ is a partition of at most $n-1$ parts.

An open problem, settled only in the $2$-variable case \cite{Evans95},
is to find (and prove) a finite field analogue of
Kadell's integral.

\medskip

The binomial theorem for Jack polynomials \eqref{BTP} 
can succinctly be written in hypergeometric notation as
\[
{_1 F_0}^{(\gamma)}\biggl(\genfrac{}{}{0pt}{}
{a}{\text{--}}\,;x\biggr)=
\prod_{i=1}^n \frac{1}{(1-x_i)^a} 
\]
where $_1 F_0^{(\gamma)}$ is an example of 
a hypergeometric function with Jack polynomial argument
\begin{equation}\label{pfq}
{_r F_s}^{(\gamma)}\biggl(\genfrac{}{}{0pt}{}
{a_1,\dots,a_r}{b_1,\dots,b_s};x\biggr) :=
\sum_{\la}
\frac{\Poch{a_1}{\la}{\gamma} \cdots \Poch{a_r}{\la}{\gamma}}
{\Poch{b_1}{\la}{\gamma} \cdots \Poch{b_s}{\la}{\gamma}}\:
\frac{P_{\la}^{(1/\gamma)}(x)}{c'_{\la}(\gamma)}.
\end{equation}
When $n=1$, so that $x$ is a scalar, this function reduces to the 
classical hypergeometric function $_r F_s$.
For general $n$, hypergeometric functions of the type \eqref{pfq}
have their genesis in the work of A.G.~Constantine \cite{Constantine63},
C.S.~Herz \cite{Herz55} and R.J.~Muirhead \cite{Muirhead70}, 
but were first studied in their full form presented
here by Kaneko \cite{Kaneko93}, A.~Kor\'anyi \cite{Koranyi}
and Yan \cite{Yan92}. The case $r=2$, $s=1$ 
of \eqref{pfq} shares many properties in common 
with its $n=1$ counterpart, the Gauss hypergeometric function.

One such property is that ${_2F_1}^{(\gamma)}$ solves the
$n$-dimensional analogue of Euler's hypergeometric equation.
Specifically, Yan \cite{Yan92} and Kaneko \cite{Kaneko93}
independently showed that
${_2 F_1}^{(\gamma)}(a,b;c;x)$ is the unique symmetric function, 
analytic at the origin, that solves the system of $n$ partial 
differential equations
\begin{multline*}
x_i (1 - x_i) \frac{\partial^2 F}{\partial x_i^2} +
\Bigl( c -  (n - 1)\gamma -
\bigl(a + b + 1 - (n-1)\gamma\bigr) x_i \Bigr)
\frac{\partial F}{\partial x_i} - abF  \\
 + \gamma\sum_{\substack{j=1 \\ j \ne i}}^n
\frac{1}{x_i - x_j} \biggl(x_i(1 - x_i)
 \frac{\partial F}{\partial x_i}
- x_j(1 - x_j)
\frac{\partial F}{\partial x_j} \biggr) = 0
\end{multline*}
for $1\leq i\leq n$.

In the one-variable theory the Gauss hypergeometric function admits
an integral representation due to Euler \cite{Euler1769}
\begin{equation}
{_2 F_1}\biggl(\genfrac{}{}{0pt}{}{a,b}{c};z\biggr)=
\frac{\Gamma(c)}{\Gamma(b)\Gamma(c-b)}
\int_0^1 t^{b-1}(1-t)^{c-b-1}(1-zt)^{-a}\,\dup t,
\end{equation}
for $\Re(c)>\Re(b)>0$ with a branch cut in the complex $z$-plane
from $1$ to infinity. When $z=1$ the integral on the right is 
the beta integral \eqref{betaInt}, resulting in a closed form
evaluation of the $_2F_1$ due to Gauss. 
In the multivariable theory an analogous results holds,
where now the key integral-evaluation is provided 
by the Selberg integral.  
Multiplying both sides of Kadell's integral by 
$z^{\abs{\la}} \Poch{a}{\la}{\gamma}/c'_{\la}(\gamma)$, 
summing the left-hand side using the binomial theorem \eqref{BTP},
and using the definition \eqref{pfq} on the right-hand side, shows that 
Euler's integral extends to \cite{Kaneko93} \label{page2F1}
\begin{multline}\label{2F1Int}
{_2 F_1}^{(\gamma)}\biggl(\genfrac{}{}{0pt}{}{a,b}{c};(z^n)\biggr) 
=\prod_{j=0}^{n-1}
\frac{\Gamma(c-j\gamma)\Gamma(1+\gamma)}
{\Gamma(b-j\gamma)\Gamma(c-b-j\gamma)\Gamma(1+(j+1)\gamma)} \\
\times
\int_0^1  \cdots \int_0^1
\prod_{i=1}^n t_i^{\alpha-1} (1-t_i)^{\beta-1}(1-z t_i)^{-a} 
\prod_{1\leq i<j\leq n} \abs{t_i-t_j}^{2\gamma} \, 
\dup t_1 \cdots \dup t_n,
\end{multline}
with $\alpha=b-(n-1)\gamma$ and $\beta=c-b-(n-1)\gamma$. 
Evaluating the $z=1$ instance of the integral by the
Selberg integral (which, incidentally, follows by taking $z=0$ in 
\eqref{2F1Int}) implies a generalized Gauss summation \cite{Yan92}
\[
{_2 F_1}^{(\gamma)}\biggl(\genfrac{}{}{0pt}{}{a,b}{c};(1^n)\biggr)
=\prod_{j=0}^{n-1}
\frac{\Gamma(c-j\gamma)\Gamma(c-a-b-j\gamma)}
{\Gamma(c-a-j\gamma)\Gamma(c-b-j\gamma)}.
\]

\medskip

In addition to relating to the Selberg integral, the Cauchy identity
also gives rise to a special limiting case of the Selberg density, 
referred to as the Laguerre PDF \cite{Johansson00,FR02}. 
To motivate the origin of this we remark that the Jack polynomial 
at $\gamma=1$ is equal to the Schur polynomial $s_{\la}$ 
while $c_{\la}(1)/c'_{\la}(1)=1$.
The normalised summand of \eqref{CaP} then reads
\[
s_{\la}(x) s_{\la}(y)\,
\prod_{i=1}^n\prod_{j=1}^m (1-x_i y_j)
\]
which may be recognised as the measure on partitions
induced by the Robinson--Schensted--Knuth correspondence,
see e.g., \cite{Okounkov01,Johansson02}.
As such the Schur measure holds a special place in certain
studies relating to the representation theory of the symmetric group
\cite{BO01,Regev81}.

For general $\gamma$ the normalised summand of \eqref{CaP} implies 
the more general measure on partitions
\begin{equation}\label{LP}
\frac{c_{\la}(\gamma)}{c'_{\la}(\gamma)}\,
P_{\la}^{(1/\gamma)}(x) P_{\la}^{(1/\gamma)}(y)\,
\prod_{i=1}^n\prod_{j=1}^m (1-x_i y_j)^{\gamma},
\end{equation}
where $n\leq m$ and $\la$ a partition of at most $n$ parts.
To obtain the Laguerre PDF \cite{Johansson00,FR02} one needs to 
specialize $x$ and $y$ in \eqref{LP} to
\[
x_i=q^{1/2},\quad 1\leq i\leq n \quad \text{and}\quad
y_j=q^{1/2},\quad 1\leq j\leq m.
\]
Use of the Jack polynomial evaluation formula \eqref{ef}
allows all terms in \eqref{LP} to be made explicit. 
The remaining step is to take the scaling limit, turning the
discrete measure on partitions into a continuous one on the
positive real line.
This is done by setting $q=\exp(-1/L)$,
introducing the scaled variables $t_j := \la_j/L$ and then 
taking the large $L$ limit for fixed $t_j$. 
One finds that \eqref{LP} multiplied by $L^n$
tends to the PDF, supported on $t_1>t_2>\cdots>t_n>0$, 
\begin{equation}\label{Ll}
\frac{1}{W_n((m-n+1)\gamma,\gamma)}
\prod_{i=1}^n \eup^{-t_i} t_i^{(m-n+1)\gamma-1}
\prod_{1 \le i < j \le n} (t_i-t_j)^{2\gamma},
\end{equation}
where $W_n$ is the normalization
\begin{equation}\label{Wn}
W_n(\alpha,\gamma)=\frac{1}{n!}\prod_{j=0}^{n-1}
\frac{\Gamma(\alpha+j\gamma)\Gamma((j+1)\gamma)}
{\Gamma(\gamma)}.
\end{equation}
To obtain this same result as a limit of the Selberg density, order
the integration variables in the latter and write $\beta=L$.
Then the change of variables $t_i \mapsto t_i/L$ followed by the
limit $L \to \infty$ results in \eqref{Ll}
after identifying $\alpha$ with $(m-n+1)\gamma$. This limiting
``Laguerre'' case of the Selberg integral, leading to the evaluation 
\eqref{Wn}, is contained as equation (1) in the letter from Selberg 
to Dyson reprinted on page \pageref{pagebrief}.

Askey and D.~Richards \cite{AR89} (see also \cite{Mehta04})
have shown that after some fairly straightforward manipulations and
a change of variables the Laguerre limit of the Selberg
integral leads to the evaluation
\begin{multline}\label{AsRi}
\int_D\; 
\prod_{i=1}^n t_i^{\alpha-1}
\Bigl(1-\sum_{i=1}^n t_i\Bigr)^{\beta-1}
\prod_{1\le i<j\le n} \abs{t_i-t_j}^{2\gamma}
\dup t_1 \cdots \dup t_n \\
=\frac{\Gamma(\beta)}{\Gamma(\alpha n+\beta+n(n-1)\gamma)}
\prod_{j=0}^{n-1} \frac{\Gamma(\alpha+j\gamma) \Gamma(1+(j+1)\gamma)}
{\Gamma(1+\gamma)},
\end{multline}
where $D$ is the domain $t_i\geq 0$ for $1\leq i\leq n$
such that $t_1+\cdots+t_n \le 1$, and 
\[
\Re(\alpha)>0,~\Re(\beta)>0,~\Re(\gamma)>
-\min\{1/n,\Re(\alpha)/(n-1)\}.
\]
According to Askey and Richards, the first statement of \eqref{AsRi} 
is due to Selberg at a meeting held in Sri Lanka during December 1987.
The intriguing point is that while Selberg did not give his derivation of 
\eqref{AsRi}, he is reported to have said that it was different to
\eqref{SelInt}, and had the advantage of working in the finite field case.
The derivation given in \cite{AR89} does not work in the finite field case,
and therefore must be different to that known to Selberg.

\subsection*{\texorpdfstring{$q$-}{}Integrals and constant terms}

Motivated by the Selberg integral and its success in dealing
with Dyson and Macdonald type constant term identities, 
Askey in 1980, was led to consider several $q$-analogues of the Selberg 
integral and to study connections to $q$-constant term identities.
In fact, one learns from \cite{Askey98} that he had earlier spent
time searching for a proof of the Mehta integral upon its appearance
in the problem section of SIAM Review.
To describe some of Askey's work we require the multiple 
Jackson or $q$-integral 
\[
\int_0^1\cdots\int_0^1 f(t)\,
\dup_q t_1\cdots \dup_q t_n := (1-q)^n
\sum_{k_1,\dots,k_n=0}^{\infty} f(q^k) q^{k_1+\cdots+k_n}
\]
with $t=(t_1,\dots,t_n)$, $q^k=(q^{k_1},\dots,q^{k_n})$ and $0<q<1$,
and where it is assumed that the multiple sum on the right is absolutely 
convergent. Also needed is the $q$-shifted factorial
\[
(a;q)_z=\prod_{j=0}^{\infty}\frac{1-aq^j}{1-aq^{z+j}}
\] 
for $z\in\Complex$. Probably the most important of the $q$-Selberg 
integrals considered by Askey is
\begin{multline}\label{JN}
S_n(\alpha,\beta,\gamma;q) \\ := 
\int_0^1 \cdots  \int_0^1 
\prod_{i=1}^n t_i^{\alpha-1} (qt_i;q)_{\beta-1}
\prod_{1 \le i < j \le n} t_i^{2\gamma}(q^{1-\gamma}t_j/t_i;q)_{2\gamma}\,
\dup_q t_1 \cdots \dup_q t_n.
\end{multline}
It is immediate, at least formally, that
\[
\lim_{q \to 1^{-}} S_n(\alpha,\beta,\gamma;q) = S_n(\alpha,\beta,\gamma).
\]
For $\Re(\alpha)>0$, $\gamma$ a nonnegative integer, say $k$, and
$\beta\in\Complex$ excluding the nonpositive integers
Askey conjectured (and proved for $k=2$) that \cite[Conjecture 1]{Askey80}
\[
S_n(\alpha,\beta,k;q)
=q^{\alpha k\binom{n}{2}+2k^2\binom{n}{3}}
\prod_{j=0}^{n-1}
\frac{\Gamma_q(\alpha+(j-1)k) \Gamma_q(\beta+(j-1)k) 
\Gamma_q(1+jk)}{\Gamma_q(\alpha+\beta+(n+j-2)k) 
\Gamma_q(1+k)},
\]
where $\Gamma_q(x)$ is the $q$-gamma function
\[
\Gamma_q(x)=\frac{(q;q)_{x-1}}{(1-q)^{x-1}}.
\]
For Askey's other $q$-Selberg integrals and many further results
relating to Jackson-integral type extensions of beta integrals,
see \cite{AA81,Aomoto95,Aomoto98,Askey80,Evans92,Ito97,
Kadell88a,Kaneko96,W05}.

In 1988 Askey's conjecture was proved independently by
L.~Habsieger \cite{Habsieger88} and K.~Kadell \cite{Kadell88}.
Both then applied the $q$-analogue of the 
identity \eqref{si} to \eqref{JN} to obtain a $q$-generalization 
of the Morris integral \eqref{Morris}. Expressing this integral 
as a constant term identity they thus proved Morris' $q$-constant
term conjecture \cite{Morris82}
\begin{multline*}
\text{CT} \, \prod_{j=1}^n
(x_j;q)_a (q/x_j;q)_b \prod_{1\le i<j\le n}
(x_i/x_j;q)_k(q x_j/x_i;q)_k \\
= \prod_{j=0}^{n-1}
\frac{\Gamma_q(1+a+b+jk) \Gamma_q (1+(j+1)k)}
{\Gamma_q (1+a+jk)\Gamma_q (1+b+jk) \Gamma_q (1+k)}.
\end{multline*}
When $a=b=0$ this is precisely the A$_{n-1}$ case
of the $q$-Macdonald constant term conjecture 
\cite[Conjecture 3.1]{Macdonald82}
\begin{subequations}\label{qMac}
\begin{equation}
\text{CT} \prod_{{\alpha} \in \Phi^{+}} \prod_{i=1}^k
(1-q^{i-1} \eup^{-\alpha})(1-q^i \eup^{\alpha}) = 
\prod_{i=1}^l \qbin{d_i k}{k}_q,
\end{equation}
where 
\[
\qbin{n}{k}_q=\prod_{i=0}^k \frac{1-q^{n-k+i}}{1-q^i}
\]
is a $q$-binomial coefficient, and $\Phi$ is a reduced (finite)
root system.
To also include the root systems of type BC one again needs
the numbers $s(\alpha)$ as defined above equation \eqref{MacCT}
\cite[Conjecture 3.4]{Macdonald82}
\begin{equation}
\text{CT} \prod_{{\alpha} \in \Phi^{+}} \prod_{i=1}^k
(1-q^{is(\alpha)-1}\eup^{-\alpha})(1-q^{(i-1)s(\alpha)+1}\eup^{\alpha})
=\prod_{i=1}^l \qbin{d_i k}{k}_q.
\end{equation}
\end{subequations}

The A$_{n-1}$ case of \eqref{qMac}, was in fact proved prior to the work
of Habsieger and Kadell by Zeilberger and D.M.~Bressoud \cite{ZB85}, who
proved the more general $q$-Dyson conjecture
\[
\text{CT} \prod_{1\leq i<j\leq n}
\Big(\frac{x_i}{x_j}\Big)_{a_i} \Big(\frac{qx_j}{x_i}\Big)_{a_j}
=\frac{(q;q)_{a_1+\cdots+a_n}}{(q;q)_{a_1}\cdots (q;q)_{a_n}}
\]
formulated by G.E.~Andrews \cite{Andrews75}.

R.A.~Gustafson \cite{Gustafson90}, at around the same time as Anderson's 
work on the Dixon--Anderson integral, made use of a further $q$-generalization 
of \eqref{GWA} and invented the same general strategy as
used in \cite{Anderson91} to derive the BC$_n$-type constant term identity 
\begin{equation}\label{CTBCn}
\text{CT} \, \Delta(x;t,t_1,\dots,t_4)
= 2^n n! \prod_{j=1}^n \frac{(t;q)_{\infty} 
(t^{n+j-2} t_1 t_2 t_3 t_4; q)_{\infty}}
{(t^j; q)_{\infty} (q;q)_{\infty} \prod_{1 \le r < s \le 4}
(t^{j-1}t_r t_s; q)_{\infty}},
\end{equation}
where
\begin{multline}\label{DeltaG}
\Delta(x;t,t_1,\dots,t_4):=
\prod_{i=1}^n \frac{(x_i^2;q)_{\infty} (x_i^{-2};q)_{\infty}}
{\prod_{r=1}^4 (t_r x_i; q)_{\infty} (t_r x_i^{-1}; q)_{\infty}} \\
\times
\prod_{1 \le i < j \le n}
\frac{ (x_i x_j^{-1}; q)_{\infty} (x_i^{-1} x_j; q)_{\infty} 
(x_i x_j;q)_{\infty} (x_i^{-1} x_j^{-1};q)_{\infty}}
{(tx_i x_j^{-1};q)_{\infty} (tx_i^{-1} x_j; q)_{\infty}
(t x_i x_j; q)_{\infty} (t x_i^{-1} x_j^{-1}; q)_{\infty}}.
\end{multline}
This is a generalization of the so-called Macdonald--Morris constant
term identity and implies the B$_n$, C$_n$, D$_n$ and BC$_n$
cases of \eqref{qMac} through specialisation \cite{Gustafson90}.

Most other cases of \eqref{qMac} were proved on a case by case basis,
often using methods based on $q$-integrals of Selberg type
\cite{Evans92,Garvan90,GG91,Habsieger86,Kadell94,Stembridge88,
Zeilberger87,Zeilberger88,Zeilberger89}, 
but the three exceptional root systems E$_6$, E$_7$ and E$_8$ 
refused to surrender until Cherednik gave a uniform proof for all 
reduced root systems based on his theory of double affine Hecke algebras 
\cite{Cherednik95}.

\medskip

Returning to Askey's $q$-Selberg integral we remark
that the density function corresponding to the integrand of 
\eqref{JN} can be deduced from Macdonald polynomial theory
\cite{Macdonald95}, following a procedure similar to that of
deducing \eqref{Ll} from \eqref{LP} \cite{FR02}. 
Macdonald polynomials $P_{\la}(x;q,t)$ are
generalizations of Jack polynomials --- the latter being reclaimed
according to
$\lim_{q \to 1} P_{\la}(x;q,q^{\gamma})=
P_{\la}^{(1/\gamma)}(x)$ ---
exhibiting the structure \eqref{stu} and the orthogonality
\[
\int_{-\pi}^{\pi}\cdots\int_{-\pi}^{\pi} P_{\la}(\eup^{i \theta};q,t)
P_{\mu}(\eup^{-i\theta};q,t) \prod_{1\leq i<j\leq n}
\frac{(\eup^{i (\theta_i - \theta_j)};q)_{\infty}}
{(t \eup^{i (\theta_i - \theta_j)};q)_{\infty}}\,
\dup \theta_1 \cdots \dup \theta_n \propto \delta_{\la\mu}.
\]
It is the connection between affine Hecke algebras and Macdonald
type orthogonal polynomials that is at the heart of Cherednik's
proof of the $q$-constant term conjectures for arbitrary root systems
\cite{Cherednik95}, see also \cite{Macdonald03}.

\subsection*{Multivariable orthogonal polynomials}
By the change of variables $t=(1-x)/2$ and a shift in
$\alpha$ and $\beta$ by $1$ the Euler
beta integral \eqref{betaInt} takes the form
\[
J(\alpha,\beta):=\int_{-1}^1 (1-x)^{\alpha} (1+x)^{\beta} \dup x =
2^{\alpha+\beta+1}\:
\frac{\Gamma(\alpha+1)\Gamma(\beta+1)}{\Gamma(\alpha+\beta+1)}.
\]
The integrand on the left is the weight function of the
Jacobi polynomials $P_n^{(\alpha,\beta)}(x)$ \cite{AAR99,Ismail05}.
Up to normalization these are the unique functions,
analytic around the origin, solving the second order ODE
\[
(1-x^2) y''+\Bigl(\beta-\alpha-x(\alpha+\beta+2)\Bigr)y'
+n(n+\alpha+\beta+1)y=0.
\]
Defining the inner product
\[
\langle f,g \rangle_{\alpha,\beta}:=
\int_{-1}^1 f(x)g(x)(1-x)^{\alpha} (1+x)^{\beta} \,\dup x
\]
the Jacobi polynomials satisfy the orthogonality relation
\[
\langle P_n^{(\alpha,\beta)},P_m^{(\alpha,\beta)}
\rangle_{\alpha,\beta}=\delta_{m,n}\:
\frac{(\alpha+1)_n(\beta+1)_n}{n! (\alpha+\beta+1)_n}\;
\frac{\alpha+\beta+1}{\alpha+\beta+2n+1}\; J(\alpha,\beta).
\]
All other classical orthogonal polynomials, such as the
Laguerre and Hermite polynomials 
(corresponding to weights $x^{\alpha}\exp(-x)$ and $\exp(-x^2)$
respectively) follow from the Jacobi polynomials
by taking appropriate limits.

Several people have studied multivariable generalizations of the 
Jacobi, Laguerre and Hermite polynomials 
\cite{BF97b,Debiard87,DX01,Lassalle91a,Lassalle91b,Lassalle91c,Macdonald,Opdam00,vanDiejen97}. 
The most general of these are the multivariable Jacobi polynomials
$P_{\la}^{(\alpha,\beta,\gamma)}(x)$
which arise as the eigenfunctions of the operator
\[
\sum_{i=1}^n \biggl( (1-x_i^2)
\frac{\partial^2}{\partial x_i^2}
+\Bigl(\beta-\alpha-x_i(\alpha+\beta+2)\Bigr)
\frac{\partial}{\partial x_i}\biggr) +
2\gamma \sum_{\substack{i,j=1 \\ i\neq j}}^n
\frac{1-x_i^2}{x_i-x_j}\,
\frac{\partial}{\partial x_i}.
\]
The $P_{\la}^{(\alpha,\beta,\gamma)}(x)$ are orthogonal with respect 
to an inner product with weight function derived from the Selberg integral.
With
\begin{equation}\label{JacobiIP}
\langle f,g \rangle_{\alpha,\beta,\gamma}:=
\Int_{[-1,1]^n} f(x)g(x)
\prod_{i=1}^n (1-x_i)^{\alpha} (1+x_i)^{\beta} 
\prod_{1\leq i<j\leq n} \abs{x_i-x_j}^{2\gamma}
\,\dup x_1 \cdots \dup x_n
\end{equation}
for $x=(x_1,\dots,x_n)$, the multivariable Jacobi polynomials
satisfy
\[
\langle P_{\la}^{(\alpha,\beta,\gamma)},P_{\mu}^{(\alpha,\beta,\gamma)}
\rangle_{\alpha,\beta,\gamma}=0 \quad \text{if $\la\neq\mu$}.
\]
The quadratic norm evaluation can be computed explicitly
in term of Pochhammer symbols and gamma functions using the
shift operators of Heckman and Opdam \cite{Opdam89}. From the
Selberg integral it of course immediately follows that
\[
\langle 1,1 \rangle_{\alpha,\beta,\gamma}=
2^{n(\alpha+\beta+1+(n-1)\gamma)}
S_n(\alpha+1,\beta+1,\gamma).
\]

Two important limiting cases of the inner product \eqref{JacobiIP}
are 
\[
\langle f,g \rangle_{\gamma}:=
\Int_{\R^n} f(x)g(x)
\prod_{i=1}^n \eup^{-x_i^2} \prod_{1\leq i<j\leq n} \abs{x_i-x_j}^{2\gamma}
\,\dup x_1 \cdots \dup x_n
\]
and
\[
\langle f,g \rangle_{\alpha,\gamma}:=
\Int_{[0,\infty)^n} f(x)g(x)
\prod_{i=1}^n x_i^{\alpha} \eup^{-x_i} 
\prod_{1\leq i<j\leq n} \abs{x_i-x_j}^{2\gamma}
\,\dup x_1 \cdots \dup x_n.
\]
The corresponding families of orthogonal
polynomials are referred to as the multivariable Hermite and Laguerre
polynomials, respectively. In particular note that
\[
\langle 1,1 \rangle_{\gamma}=
2^{\binom{n}{2}\gamma} \pi^{n/2} \, F_n(\gamma)
\]
with $F_n$ the Mehta integral \eqref{FNz}, and
\[
\langle 1,1 \rangle_{\alpha,\gamma}=
n! \, W_n(\alpha+1,\gamma),
\]
with $W_n$ given by \eqref{Wn}.

All of the orthogonal polynomials mentioned above admit 
$q$-analogues. In the $q$-theory the role of the Selberg integral
is played by Askey's $q$-Selberg integral \eqref{JN}
and generalisations thereof.
In the case of the Jacobi
polynomials these $q$-analogues are known as the multivariable
big and little $q$-Jacobi polynomials, and were introduced
by J.V.~Stokman \cite{Stokman97}.
Stokman \cite{Stokman00} also showed how the big and little
$q$-Jacobi polynomials arise as special
limits of the Koornwinder polynomials \cite{Koornwinder92}.
The latter are multivariable analogues of the Askey--Wilson orthogonal
polynomials \cite{AW85} and may be viewed as the generalizations 
of the Macdonald polynomials to the root system BC$_n$.
The relevant inner product in this case is given by 
\cite{Koornwinder92,vanDiejen96}
\begin{equation}\label{KIP}
\langle f,g\rangle_{t,t_1,\dots,t_4}:=
\frac{1}{(2\pi)^n}
\int_{-\pi}^{\pi}\cdots
\int_{-\pi}^{\pi}f(\eup^{\iup x})g(\eup^{\iup x})
\Delta(\eup^{\iup x};t,t_1,\dots,t_4)
\,\dup x_1\cdots\dup x_n,
\end{equation}
where $\exp(\iup x)=(\exp(\iup x_1),\dots,\exp(\iup x_n))$,
$\Delta(x;t,t_1,\dots,t_4)$ is the weight function \eqref{DeltaG}
of Gustafson's constant term identity, and
$f$ and $g$ are BC$_n$ symmetric functions.
($f(z)$ is BC$_n$ symmetric if $f(\exp(\iup x))$ is symmetric
under signed permutations of $x=(x_1,\dots,x_n)$.)
The evaluation of $\langle 1,1\rangle_{t,\dots,t_4}$
is of course provided by the right-hand side of \eqref{CTBCn}.

\section{Recent and current research directions}

\subsection*{The case of \texorpdfstring{$\gamma$}{} a positive integer}

Two recent studies have identified special features of
the Selberg integral for $\gamma$ a positive integer. The first of these
is due to J.-G.~Luque and J.-Y.~Thibon \cite{LT03}, and exhibits an 
inter-relation with a special class of hyperdeterminants. 
The second, due to Stanley \cite{Stanley05},
gives a probabilistic interpretation of Selberg's integral.

For a $k$th order tensor $A = [A_{i_1 i_2 \cdots i_k}]$ on an $n$-dimensional
space (so that $1 \le i_p \le n$) the hyperdeterminant
has been defined by Cayley (see references in \cite{LT03}) as
\[
{\det}_k(A)
:=\sum_{\sigma_2,\dots,\sigma_k\in\Symm_n}
\epsilon(\sigma_2)\cdots\epsilon(\sigma_k)
\prod_{i=1}^n A_{i,\sigma_2(i),\dots,\sigma_k(i)}
\]
where $\epsilon(\sigma)$ denotes the signature of the permutation $\sigma$.
For $k$ odd this vanishes while $k=2$ corresponds to the usual definition
of a determinant. 

For an arbitrary measure $\mu(x)$ it is easy to see by use of the 
Vandermonde determinant formula that in the so-called Hankel case
\[
A_{i_1 i_2 \cdots i_k} = \mu_{i_1 + i_2 + \cdots + i_k-k}, \qquad
\mu_j := \int_{-\infty}^\infty x^j \, \dup \mu(x)
\]
the corresponding hyperdeterminant is equal to a multiple integral,
\[
{\det}_{2k}(A)= \frac{1}{n!}
\int \cdots \int \prod_{1 \le i < j \le n}
(x_i-x_j)^{2k} \, \dup \mu(x_1) \cdots \dup \mu(x_n).
\]
For $\dup\mu(x) = x^{\alpha - 1} (1 - x)^{\beta - 1}\dup x$ 
on $x\in[0,1]$ this is precisely the Selberg integral
with $\gamma=k$ an integer.

\medskip

The probabilistic interpretation of the Selberg integral given in
\cite{Stanley05} applies for  $\alpha = \beta = 1$ and 
$\gamma$ a positive integer.
In fact, as communicated to us by Stanley \cite{Stanley07},
this same probabilistic interpretation, extended to
$\alpha, \beta$ general nonnegative integers, is already implied by
appropriately interpreting supplementary problem I.25 of
\cite{Stanley97}.
Following \cite{Stanley07}, the setting is to choose labelled points
independently and with uniform distribution from the interval
$[0,1]$. Specifically, for each $1 \le p \le n$ and $(i,j)$ such
that $1 \le i < j \le n$, choose $\alpha - 1$ points labelled $y_p$,
$\beta - 1$ points labelled $z_p$, $n$ points labelled $t$ and $2 \gamma$
points labelled $a_{ij}$. Let $t_i$ be the $i$-th smallest point labelled
$t$. The probability that any one of the points labelled $y_p$ is to the left
of $t_i$ is $1 - t_i$; the probability that any one of the points labelled
$a_{ij}$ is between $t_i$ and $t_j$ for $i < j$ is $t_i - t_j$. It follows
immediately that the Selberg integral $S_n(\alpha,\beta,\gamma)$ is the
probability, $P_S$ say, that all the points labelled $y_p$ are to the left of
$t_p$, and all the points labelled $z_p$ are to the right of $t_p$, and all
the points labelled $a_{ij}$ lie between $t_i$ and $t_j$. 
Note that this statement remains valid for $2\gamma$ an odd integer.
An equivalent formulation (the one given in \cite{Stanley97}) is to regard the
selection of the labelled points uniformly and independently from $[0,1]$
as a random re-arrangement of the symbols themselves
(i.e.~the $t$'s, $y_p$'s, $z_p$'s and $a_{ij}$'s), in which case $P_S$
corresponds to the probability that the re-arrangement complies with
the prescribed rule.

\subsection*{Random matrix theory}

A number of interplays between random matrix theory and the Selberg integral
appearing in papers published in the last few years were discussed previously
under the heading of the Dixon--Anderson integral. 
Below two further applications of the Selberg integral to random matrix 
theory as they have occurred in current works will be outlined.

The first of these is a study by Forrester and Rains \cite{FR07}
focusing on the family of multi-dimensional integrals
\begin{multline}\label{pal}
I_{n,p}(x) \\ :=
\int_{[0,x]^p} \int_{[x,1]^{n-p}} 
\prod_{i=1}^n t_i^{\alpha-1} (1-t_i)^{\beta-1}
\abs{x-t_i}^{\tau-1} \prod_{1\leq i<j\leq n}
\abs{t_i-t_j}^{2\gamma} \, \dup t_1 \cdots \dup t_n
\end{multline}
for $0\leq p\leq n$ and $x\in[0,1]$.
In the case $\tau=1$ the integral \eqref{pal} is proportional to the
probability that for the Selberg density, interpreted as an eigenvalue PDF,
there are $p$ eigenvalues in the interval $[0,x]$ and $n-p$ eigenvalues 
in the interval $[x,1]$; 
for $\tau =1+2\gamma$ the integral \eqref{pal} relates to the 
derivative of this quantity.  The case $\tau = 2\gamma = 1$ of this was first
studied in the mathematical statistics literature for its relevance
to canonical correlation analysis \cite{Davis72}.

Theory connecting $I_{n,p}$ to a certain Fuchsian differential 
equation \cite{Davis72,Fo93,Mi93} implies that the integral is expressible 
as a linear combination of Frobenius solutions. These
are solutions to the differential equation of the form
\[
g_i(x) = x^{\sigma_i} \sum_{k=0}^\infty a_{i,k} \, x^k, \qquad
\sigma_i = i(\alpha - 1 + \tau + (i-1)\gamma)
\]
for $0\leq i\leq n$, and are normalized such that
\begin{equation}\label{rN}
\lim_{x \to 0} \frac{g_i(x)}{I_{n,p}(x)} = 1 \quad
\text{for}\quad \Re(\sigma_i) > 0. 
\end{equation}
A basic task --- essentially equivalent to finding the monodromy matrix
for the basis of integral solutions of the matrix Fuchsian system, 
of which \eqref{pal} forms the top row ---
is to give the explicit form of the coefficients $c_{p,i}$ in the expansion
\begin{equation}
I_{n,p}(x) = \sum_{i=0}^n c_{p,i} \, g_i(x).
\end{equation}

One approach to this problem is to seek a regime in parameter space 
such that for $x \to 0$ the leading behaviour of $I_{n,p}(x)$ is 
proportional to $x^{\sigma_i}$. This is achieved by changing 
variables $t_j = x u_j$ for $1\leq j\leq i$ where $p+1\leq i\leq n$.
A simple scaling of the integrand then shows that
\[
I_{n,p}(x) \sim x^{\sigma_i} S_{n-i}(\alpha+\tau+2\gamma-1,\beta,\gamma)
S_{i,p}(\alpha,\tau,\gamma),
\]
where $S_{n-i}$ is the Selberg integral and $S_{i,p}$ the 
Dotsenko--Fateev integral \eqref{DF}.
Recalling the normalization
\eqref{rN} and the recurrence \eqref{DF1} allows the sought 
coefficients to be calculated as
\begin{align*}
c_{p,i}&=0, && 0\leq i\leq p-1 \\
c_{p,i}&=(-1)^{i-p} \prod_{j=1}^{i-p}
\frac{\sin\pi(i-j+1)\gamma \sin\pi(\alpha+(i-j)\gamma)}
{\sin\pi j\gamma \sin\pi(\alpha+\tau-1+(2i-j-1)\gamma)},
&& p\leq i\leq n
\end{align*}
thus solving the problem at hand.

\medskip

The second of the applications stems from a question posed by 
B\'alint Vir\'ag at a recent AMS--IMS--SIAM summer 
research conference \cite{AIS07}. 
By way of background to his question, let us
recall a result of Mehta and Dyson \cite{MD63} which gives the 
circular ensembles identity $\alt(\text{COE}_{2n})
= \text{CSE}_n$. Here $\alt$ is the operation of integrating out 
every second eigenvalue, and the subscripts on the names of the ensembles
indicate the total number of eigenvalues. 
Let us more generally introduce the notation
CE$_{\beta,n}$ for the PDF \eqref{Cp}. 
The question posed by Vir\'ag was to investigate
extensions of the result of Mehta and Dyson, 
in which blocks of eigenvalues in CE$_{\beta,n}$ are integrated out 
to obtain another circular ensemble CE$_{\beta',n'}$. 
The Selberg integral is relevant for this purpose. 

Let $p(k;s;\text{CE}_{\beta,n})$ denote the PDF for the spacing
between eigenvalues which are $(k+1)$-st neighbours in the ensemble
CE$_{\beta,N}$. Let $\alt_m(\text{CE}_{\beta,n})$ 
denote the joint marginal distribution of every
$m$-th eigenvalue in CE$_{\beta,n}$. With this notation, if it were true that
\begin{equation}\label{pus}
\alt_m(\text{CE}_{\beta,mn}) = \text{CE}_{\beta',n}
\end{equation}
for some $m,\beta,\beta'$, then 
\begin{equation}\label{ps}
p(mk+m-1;s;\text{CE}_{\beta,mn})=
p(k;s;\text{CE}_{\beta',n}).
\end{equation}

Now the $k$-point correlation $\rho_k$ is obtained from \eqref{Cp} by 
integrating out the variables 
$\theta_{k+1},\dots,\theta_n$ and multiplying by $n!/(n-k)!$. 
It follows from this definition that $\rho_k$ is related to 
the small $s$ expansion of $p$ according to
\begin{equation}\label{ps2}
p(k;s;\text{CE}_{\beta,n}) \sim \frac{2\pi}{nk!}
\int_0^s \cdots \int_0^s 
\rho_{k+2}(0,s,\theta_1,\dots,\theta_k)\,
\dup \theta_1 \cdots \dup \theta_k.
\end{equation}
But for $\theta_1,\dots,
\theta_n$ small  the definition of $\rho_k$ also implies that
\begin{equation}\label{ps1}
\rho_k(\theta_1,\dots,\theta_k) \sim 
\frac{1}{(2\pi)^k}\, \frac{n!}{(n-k)!}\,
\frac{M_{n-k}(k \beta/2, k\beta/2,\beta/2)}{M_n(0,0,\beta/2)}
\prod_{1\le i<j\le k} \abs{\theta_i-\theta_j}^{\beta},
\end{equation}
where $M_n$ refers to the Morris integral \eqref{Morris}.
Substituting \eqref{ps1} in \eqref{ps2} and scaling the integrand 
an example of the Selberg integral is obtained, giving the formula
\begin{multline}\label{ps3}
p(k;s;\text{CE}_{\beta,n}) \sim 
\frac{1}{(2\pi)^{k+1}} \, \frac{(n-1)!}{k!(n-k-2)!} \:
s^{k+\beta(k+2)(k+1)/2} \\ 
\times
\frac{M_{n-k-2}((k+2) \beta/2, (k+2)\beta/2,\beta/2)}{M_n(0,0,\beta/2)} \:
S_{k}(\beta+1,\beta+1,\beta/2).
\end{multline}

Using the gamma function evaluations \eqref{SelInt} and \eqref{MorrisEval}, 
together with the duplication formula for the gamma function, 
one can check that in the case $m=r+1$, $\beta=2/(r+1)$ and 
$\beta'=2(r+1)$, \eqref{ps} is compatible with \eqref{ps3}. 
Thus, this investigation
based on the Selberg integral reveals parameters for which the validity of
\eqref{pus} may be expected. One can in fact proceed further and prove,
using a generalization of the Dixon--Anderson integral,
that for these parameters \eqref{pus} is indeed valid \cite{Fo07}.

\subsection*{KZ equations and the Mukhin--Varchenko conjecture}
On page \pageref{page2F1} we have seen that hypergeometric integrals
of Selberg type arise naturally as solutions of (systems) of partial
differential equations.
There is a well-developed theory extending much of this to the
setting of partial differential equations --- referred to as
Knizhnik--Zamolodchikov (KZ) equations --- based on Lie algebras
\cite{EFK03,SV91,MV00,Varchenko03}.

Let $\g$ be a simple Lie algebra of rank $n$, with simple roots, fundamental
weights and Chevalley generators given by $\ba_i$, 
$\Lambda_i$ and $e_i,f_i,h_i$ for $1\leq i\leq n$.
Let $V_{\la}$ and $V_{\mu}$ be highest weight representations of $\g$
with highest weights $\la$ and $\mu$, and let $u=u(z,w)$ be 
a function taking values in $V_{\la}\otimes V_{\mu}$ solving the
KZ equation
\[
\kappa\,\frac{\partial u}{\partial z}=\frac{\Omega}{z-w}\, u,
\qquad\quad
\kappa\,\frac{\partial u}{\partial w}=\frac{\Omega}{w-z}\, u,
\]
where $\Omega\in \g\otimes \g$ is the Casimir element.
(For the sake of simplicity
we only consider the KZ equation in two variables, $z$ and $w$; for
the more general case of $p$ variables $z_1,\dots,z_p$, see e.g.,
\cite{Varchenko03}.)
Let $\text{Sing}_{\la,\mu}[\nu]$ denote the space of singular 
vectors of weight $\nu$ in $V_{\la}\otimes V_{\mu}$
\[
\text{Sing}_{\la,\mu}[\nu]:=\{v\in V_{\la}\otimes V_{\mu}: h_i v=\nu(h_i) v,~
e_i v=0,~1\leq i\leq n\}.
\]
Then, according to a theorem of V.V.~Schechtman and Varchenko \cite{SV91},
solutions $u$ with values in 
$\text{Sing}_{\la,\mu}[\la+\mu-\sum_{i=1}^n k_i \ba_i]$ 
are expressible in terms of multiple hypergeometric integrals
\[
u(z,w)=\sum u_{\IJ}(z,w) \, f^I v_{\la}\otimes f^J v_{\mu}
\]
with coordinate functions $u_{IJ}$ given by
\[
u_{\IJ}(z,w)=\int_{\gamma}\Phi^{1/\kappa}(z,w;t)A_{\IJ}(z,w;t)\,
\dup t_1\cdots \dup t_k.
\]
Here $k:=k_1+\cdots+k_n$, $t:=(t_1,\dots,t_k)$, 
the sum is over all ordered multisets $I$ and $J$ with elements
taken from $\{1,\dots,n\}$ such that their union contains the
number $i$ exactly $k_i$ times, $v_{\la}$ and $v_{\mu}$ are the 
highest weight vectors of $V_{\la}$ and $V_{\mu}$,
$f^{I} v:=(\prod_{i\in I} f_i) v$ and $\gamma$ is a suitable 
integration domain.

The functions $\Phi$ and $A_{\IJ}$ in the integrand of $u_{\IJ}$ are
explicitly known. $A_{\IJ}$  is a rational function whose general form
is too involved to explicitly state here (an example will be given below),
and the function $\Phi$, known as the master function, is
defined as follows. The first $k_1$ integration variables are
attached to the simple root $\ba_1$, the next $k_2$ integration variables
are attached to the simple root $\ba_2$, and so on, such that
$\ba_{t_j}:=\ba_i$ if $k_1+\cdots+k_{i-1}<j\leq k_1+\cdots+k_i$. 
With this understood
\begin{multline*}
\Phi(z,w;t)=(z-w)^{\bil{\la}{\mu}}
\prod_{i=1}^k (t_i-z)^{-\bil{\la}{\ba_{t_i}}}
(t_i-w)^{-\bil{\mu}{\ba_{t_i}}} \\
\times
\prod_{1\leq i<j\leq k}(t_i-t_j)^{\bil{\ba_{t_i}}{\ba_{t_j}}},
\end{multline*}
where $\bil{\,}{}$ is the bilinear symmetric form on $\h^{\ast}$
(the dual of the Cartan subalgebra $\h$)
normalised such that $\bil{\theta}{\theta}=2$ for the
maximal root $\theta$.

The simplest possible example of a KZ solution contained 
in the Schechtman--Varchenko theorem corresponds to the rank 
$1$ Lie algebra $\g=\gsl_2=\text{A}_1$,
with simple root $\alpha_1$ and fundamental weight $\Lambda_1=\alpha_1/2$.
Taking $\la=m_1 \Lambda_1$ and $\mu=m_2 \Lambda_1$ it follows that
$u(z,w)$ takes values in the space of singular vectors 
of weight $(m_1+m_2-2k_1)\Lambda_1$. Since $n=1$ it follows that
$I=\{1^r\}$ and $J=\{1^{k_1-r}\}$ with $1\leq r\leq k_1$, so that
$u_{\IJ}$, $A_{\IJ}$ and $f^I$ can simply be denoted by 
$u_r$, $A_r$ and $f^r$. (In the case of rank one there is
no need for the index in $f_1$).
Using $n$ instead of $k_1$ (so that $n$ no longer denotes 
the rank of the Lie algebra) and writing $v_1$ and $v_2$ 
instead of $v_{\la}=v_{m_1\Lambda_1}$ and $v_{\mu}=v_{m_2\Lambda_2}$
one finds, upon the assumption that $z<w$ are both real,
\[
u(z,w)=\sum_{r=0}^n u_r(z,w) \, f^r v_1\otimes f^{n-r} v_2
\]
with
\begin{multline*}
u_r(z,w)=(z-w)^{m_1m_2/(2\kappa)}
\int_{\gamma}A_r(z,w;t) \prod_{i=1}^n (t_i-z)^{-m_1/\kappa}
(t_i-w)^{-m_2/\kappa} \\ \times
\prod_{1\leq i<j\leq n}(t_i-t_j)^{2/\kappa}\,
\dup t_1\cdots \dup t_n.
\end{multline*}
Here the domain of integration is the simplex 
$\gamma=\{t\in\R^n|~z\leq t_n\leq \cdots\leq t_1\leq w\}$,
and the rational function $A_r(z,w;t)$ is given by
\[
A_r(z,w;t)=\sum_{\substack{I\subseteq \{1,\dots,n\} \\ \abs{I}=r}}
\prod_{i\in I}\frac{1}{t_i-z}
\prod_{i\not\in I}\frac{1}{t_i-w}.
\]
The coordinate functions $u_r$ are easily recognised as 
generalizations of the Selberg integral.
In fact, for the extremal cases $r=0$ and $r=n$ they are exactly the
Selberg integral. When $r=0$, for instance,
\begin{multline*}
u_0(z,w)=(z-w)^{m_1m_2/(2\kappa)}
\int_{\gamma}\prod_{i=1}^n (t_i-z)^{-m_1/\kappa}
(t_i-w)^{-m_2/\kappa-1} \\ \times
\prod_{1\leq i<j\leq n}(t_i-t_j)^{2/\kappa}\,
\dup t_1\cdots \dup t_n.
\end{multline*}
Making the change of variables $t_i=(w-z)s_i+z$ for $1\leq i\leq n$
this yields
\[
u_0(z,w)=\frac{(-1)^A (z-w)^B}{n!} \, 
S_n\Bigl(1-\frac{m_1}{\kappa},-\frac{m_2}{\kappa},
\frac{1}{\kappa}\Bigr),
\]
where $A=n(n-1-m_1)/\kappa+n$ and 
$B=(m_1 m_2-2n(m_1+m_2)+2n(n - 1))/(2\kappa)$.

\medskip

In 2000 E.~Mukhin and Varchenko \cite{MV00} formulated a surprising
conjecture regarding the scaled master function
\[
\Phi(t)=
\prod_{i=1}^k t_i^{-\bil{\la}{\ba_{t_i}}}
(1-t_i)^{-\bil{\mu}{\ba_{t_i}}}
\prod_{1\leq i<j\leq n}(t_i-t_j)^{\bil{\ba_{t_i}}{\ba_{t_j}}}.
\]
They conjectured that if the space
$\text{Sing}[\la+\mu-\sum_{i=1}^n k_i \ba_i]$ 
of singular vectors is one-dimensional, then
\begin{equation}\label{MV}
\int \abs{\Phi(t)}^{1/\kappa} \, \dup t_1\cdots\dup t_k
\end{equation}
is expressible as a product of gamma functions.
Neither the exact integration domain nor the specific 
form for the product of gamma functions is contained in
the Mukhin--Varchenko conjecture.

For $\g=\gsl_2=\text{A}_1$ the conjecture corresponds to the
evaluation of the Selberg integral.
For $\g=\gsl_{n+1}=\text{A}_n$, and 
$\text{Sing}_{\la,\mu}[\la+\mu-\sum_{i=1}^n \ba_i]$
with $\la=\Lambda_1$, $\mu=\sum_{i=1}^n \mu_i\Lambda_i$
the conjecture simply follows by iterating the beta integral
\eqref{betaInt}, see \cite{MV00}.  
For $\g=\text{B}_n,\:\text{C}_n$ or D$_n$ and
\[
\text{Sing}_{\Lambda_1,\Lambda_1}\Bigl[2\Lambda_1-r\ba_{n-1}-s\ba_n-
\sum_{i=1}^{n-2} \ba_i\Bigr]\quad
\text{with} \quad
(r,s)=
\begin{cases}
(2,2) &\text{B}_n \\
(2,1) &\text{C}_n \\
(1,1) &\text{D}_n
\end{cases}
\]
(corresponding to the tensor product of the vector representation
of $\g$) Mimachi and T.~Takamuki \cite{MT05} established the Mukhin--Varchenko
conjecture iterating the Selberg integral for $n=2$ (B$_n$ case)
or the beta integral (C$_n$ and D$_n$ cases).

In 2003 V.~Tarasov and Varchenko employed KZ equations and the closely
related dynamical equations to settle the conjecture for
$\g=\gsl_3=\text{A}_2$. In recent work by Warnaar \cite{W07b,W07c} an 
approach to the $\gsl_{n+1}=\text{A}_n$ 
case of the Mukhin--Varchenko conjecture
was developed, based on the theory of Macdonald polynomials
and generalized hypergeometric series.
Specifically, the integral \eqref{MV} for $\g=\text{A}_n$
can be evaluated in closed form when
$\la=\la_n\Lambda_n$ and $\mu=\sum_i \mu_i\Lambda_i$
(or when $\la=\la_1\Lambda_1$ and $\mu=\sum_i \mu_i\Lambda_i$).
Stripping the integral from its Lie algebra notation
and using $\alpha_i$ and $\beta_i$ ($1\leq i\leq n$)
for exponents in the integral (so that the $\alpha_i$ no 
longer denote the simple roots)
the $\text{A}_n$ Selberg integral can be stated explicitly
as
\begin{align*}
&\int
\prod_{s=1}^n 
\biggl[ \Abs{\Delta\bigl(t^{(s)}\bigr)}^{2\gamma}
\prod_{i=1}^{k_s} \bigl(t_i^{(s)}\bigr)^{\alpha_s-1}
\bigl(1-t_i^{(s)}\bigr)^{\beta_s-1}\biggr] 
\prod_{s=1}^{n-1}
\Abs{\Delta\bigl(t^{(s)},t^{(s+1)}\bigr)}^{-\gamma} 
\; \dup t  \\ 
&=\prod_{1\leq s\leq r\leq n} \prod_{i=1}^{k_s-k_{s-1}}
\frac{\Gamma(\beta_s+\cdots+\beta_r+(i+s-r-1)\gamma)}
{\Gamma(\alpha_r
+\beta_s+\cdots+\beta_r+(i+s-r+k_r-k_{r+1}-2)\gamma)} \\
&\quad\times
\prod_{s=1}^n \prod_{i=1}^{k_s}
\frac{\Gamma(\alpha_s+(i-k_{s+1}-1)\gamma)
\Gamma(i\gamma)}{\Gamma(\gamma)}.
\end{align*}
Here $k_1,\dots,k_{n+1}$ are nonnegative integers such that
$k_{n+1}=0$ and $k_1\leq k_2\leq \dots\leq k_n$,  the exponents
$\alpha_1,\dots,\alpha_n,\beta_1,\dots,\beta_n,\gamma\in\Complex$ such 
that $\alpha_1=\cdots=\alpha_{n-1}=1$ and such that both sides 
of the identity are well-defined. 
Furthermore, $t^{(s)}=(t_1^{(s)},\dots,t_{k_s}^{(s)})$
is the set of variables attached to the $s$th simple root of $\text{A}_n$,
\[
\Delta(u)=\prod_{1\leq i<j\leq l_u}(u_i-u_j) \qquad\text{and}\qquad
\Delta(u,v)=\prod_{i=1}^{l_u}\prod_{j=1}^{l_v} (u_i-v_j)
\]
for sets of variables $u=(u_1,\dots,u_{l_u})$ and 
$v=(v_1,\dots,v_{l_v})$,
and $\dup t=\dup t^{(1)} \cdots \dup t^{(n)}$ with
$\dup t^{(s)}=\dup t_1^{(s)}\cdots \dup t_{k_s}^{(s)}$
so that the integral is $(k_1+\cdots+k_n)$-dimensional.

Not yet specified in the $\text{A}_n$ Selberg integral is 
the domain of integration, which, unfortunately, is rather
involved. A key ingredient is the set of maps 
\[
M_s:\{1,\dots,k_s\}\to \{1,\dots,k_{s+1}\}
\]
such that 
\[
M_s(i)\leq M_s(i+1) \quad\text{and}\quad
1\leq M_s(i)\leq k_{s+1}-k_s+i.
\]
A standard counting argument shows that there are exactly
$c_{k_{s+1},k_s}$ admissible $M_s$, where $c_{n,k}$ is the
row $(n,k)$ entry in the Catalan triangle, or, 
equivalently, the number of standard Young tableaux of shape $(n,k)$.
Given $M_s$ fix an ordering among the $t_i^{(s)}$ and $t_j^{(s+1)}$
as 
\begin{equation}\label{ineq}
t^{(s+1)}_{M_s(i)}\leq t_i^{(s)}\leq t^{(s+1)}_{M_s(i)-1}
\quad \text{for~~$1\leq i\leq k_s$},
\end{equation}
where $t^{(s+1)}_0:=\infty$.
Given admissible maps $M_1,\dots,M_{n-1}$ define
$D_{M_1,\dots,M_{n-1}}^{k_1,\dots,k_n}$
as the set of points
\[
(t^{(1)}_1,\dots,t^{(1)}_{k_1},t^{(2)}_1,\dots,t^{(2)}_{k_2},
\dots,t^{(n)}_1,\dots,t^{(n)}_{k_n})
\]
such that \eqref{ineq} holds for all $1\leq s\leq n-1$ and 
\begin{equation*}
0\leq t^{(s)}_{k_s}\leq \dots\leq t^{(s)}_1\leq 1
\end{equation*}
holds for all $1\leq s\leq n$.
Then the domain of integration, written as a chain, is given by
\[
\sum_{M_1,\dots,M_{n-1}}
F_{M_1,\dots,M_{n-1}}^{k_1,\dots,k_n}(\gamma)
D_{M_1,\dots,M_{n-1}}^{k_1,\dots,k_n},
\]
where
\begin{equation*}
F_{M_1,\dots,M_{n-1}}^{k_1,\dots,k_n}(\gamma)=
\prod_{s=1}^{n-1}\prod_{i=1}^{k_s}
\frac{\sin\bigl(\pi(i+k_{s+1}-k_s-M_s(i)+1)\gamma\bigr)}
{\sin\bigl(\pi(i+k_{s+1}-k_s)\gamma\bigr)}.
\end{equation*}

In complete analogy with the ordinary Selberg integral, the evaluation
of the A$_n$ Selberg integral can be generalized to include a Jack
polynomial in the integrand, thus generalizing the Kadell integral
\eqref{KadellInt}, see \cite{W07c}.

\subsection*{Elliptic Selberg integrals} 

In the last few years there has been rapid progress in the field of 
elliptic generalizations of hypergeometric series, see \cite{GR04,Spiridonov}.
Classical hypergeometric series $\sum_{n=0}^\infty c_n$ are 
characterized by the ratio $c_{n+1}/c_n$ being a
rational function of $n$. Their elliptic counterparts have the same 
ratio equal to an elliptic function of $n$.

As well as the classical hypergeometric series permitting elliptic
generalizations, so do related integrals such as the Euler beta integral
\eqref{betaInt}. In the elliptic theory the ordinary gamma function
must be replaced by what is known as the elliptic gamma function
\begin{equation}\label{ellG}
\Gamma(z;p,q)=\prod_{i,j=0}^{\infty}
\frac{1-z^{-1}p^{i+1}q^{j+1}}{1-zp^iq^j},
\end{equation}
defined for $\abs{p},\abs{q}<1$. This function can be traced back 
to E.W.~Barnes in 1904 \cite{Barnes04}, but was given prominence
through the recent work of S.N.M.~Ruijsenaars \cite{Ruijsenaars97}. 
It permits the extension of the standard gamma recurrence to
\[
\Gamma(qz;p,q)=\theta(z;p)\Gamma(z;p,q),
\]
where $\theta(z;p)=(z;p)_{\infty}(p/z;p)_{\infty}$ is a normalised
theta function. Another fundamental property of the elliptic gamma 
function is the functional equation
\begin{equation}\label{func}
\Gamma(z;p,q)=\frac{1}{\Gamma(pq/z;p,q)} 
\end{equation}
which follows immediately from the definition \eqref{ellG}.

The elliptic analogue of the beta integral \eqref{betaInt}
was discovered in 2000 by V.P.~Spi\-ri\-do\-nov \cite{Spiridonov01}
\begin{equation}\label{Sr}
\int_{C} \frac{\prod_{r=1}^6 
\Gamma(t_r z^{\pm 1}; p,q)}{\Gamma(z^{\pm 2};p,q)}
\, \frac{\dup z}{2\pi\iup z} =
\frac{2}{(p;p)_{\infty} (q;q)_{\infty}}\,
\prod_{1 \le r < s \le 6} \Gamma(t_r t_s;p,q),
\end{equation}
where each $\abs{t_r} < 1$, $C$ is the positively oriented
unit circle, $\prod_{r=1}^6 t_r=pq$ and
\[
\Gamma(tz^{\pm m}; p,q) := \Gamma(tz^m;p,q) \Gamma(tz^{-m};p,q).
\]
The $p\to 0$ limit is the well-known Rahman \cite{NR85,Rahman86} integral,
which itself is an extension of the Askey--Wilson integral \cite{AW85}.
For the reduction of this last integral to the beta integral \eqref{betaInt}
see \cite{GR04}.

J.F.~van Diejen and Spiridonov \cite{DS01} have given an $n$-dimensional
generalization of \eqref{Sr} which may be viewed as an elliptic
extension of the Selberg integral. This integral, the $p\to 0$ limit
of which was first obtained by Gustafson \cite{Gustafson94}, takes the form
\begin{multline}\label{Sp2}
\int_{C^n}
\prod_{1 \le i < j \le n}
\frac{\Gamma(t z_i^{\pm 1} z_j^{\pm 1}; p,q)}
{\Gamma(z_i^{\pm 1} z_j^{\pm 1}; p,q)}
\prod_{i=1}^n \frac{\prod_{r=1}^6 \Gamma(t_r  z_i^{\pm 1}; p,q)}
{\Gamma(z_i^{\pm 2}; p,q)} \,
\frac{\dup z_1}{2 \pi \iup z_1} \cdots \frac{\dup z_n}{2 \pi \iup z_n} \\
=\frac{2^n n!}{(p;p)^n_{\infty} (q;q)^n_{\infty}}\,
\prod_{j=1}^n \biggl(
\frac{\Gamma(t^j;p,q)}{\Gamma(t;p,q)}
\prod_{1 \le r < s \le 6}
\Gamma(t^{j-1} t_r t_s;p,q) \biggr),
\end{multline}
where $\abs{t},\abs{t_1},\dots,\abs{t_6}<1$ and
$t^{2n-2} \prod_{r=1}^6 t_j = pq$.
van Diejen and Spiridonov provided a proof of \eqref{Sp2} along 
the lines of the Anderson and Gustafson proofs of \eqref{SelInt} and
\eqref{CTBCn} respectively. This required an elliptic 
generalization of the Dixon--Anderson integral \eqref{GWA} which,
initially, was proved making an assumption about the
vanishing of certain elliptic integrals. 
A complete proof of the elliptic Dixon--Anderson integral was
subsequently given by Rains \cite{Rains07} and by Spiridonov
\cite{Spiridonov07}.

The reduction of the elliptic Selberg integral \eqref{Sp2} to the
ordinary Selberg integral is rather cumbersome, requiring several limits,
variables changes and specializations of parameters
\cite{Rains06}. Fairly straightforward,
however, is to see that \eqref{Sp2} provides an elliptic extension of 
Gustafson's BC$_n$ constant term identity \eqref{CTBCn}. To see this
one first needs to eliminate $t_6$ using $t^{2n-2} \prod_{r=1}^6 t_r = pq$.
This gives rise to several elliptic gamma functions of the form $\Gamma(pq A)$
which, by \eqref{func}, may be replaced by $1/\Gamma(1/A)$. 
After these elementary manipulations the $p\to 0$ limit can be carried out,
using that $\Gamma(z;0,q)=1/(z;q)_{\infty}$. Finally taking $t_5=0$ and
interpreting the resulting integral as a constant term identity yields
\eqref{CTBCn}.

Analogous to \eqref{KIP}, the integrand of \eqref{Sp2}
can be used to define an inner product. Rains \cite{Rains07}
has specified a family of abelian functions which are biorthogonal
with respect to this inner product, extending the Rahman--Spiridonov
theory \cite{Rahman86,Spiridonov01b} of such functions to 
the multivariable setting, as well as generalising the Koornwinder 
polynomials and Okounkov BC$_n$ interpolation polynomials \cite{Okounkov98b}
to the elliptic level. (These functions were independently
introduced by H.~Coskun and Gustafson in \cite{CG06} without the use
of elliptic Selberg type integrals.) 
Rains also extended the integrand of \eqref{Sp2} analogous to the 
$_2F_1$ extension \eqref{2F1Int} of the Selberg integral,
and obtained transformation formulas for the resulting 
elliptic hypergeometric integrals.
By considering the reduction of his theory to the Selberg level
Rains obtained, for example, \cite{Rains}
\begin{multline*}
\int^1_0 \cdots \int^1_0 
P_{\la}^{(1/\gamma)}(t) P_{\mu}^{(1/\gamma)}(t)
\prod_{i=1}^n t_i^{\alpha-1} (1-t_i)^{\gamma-1}
\prod_{1\leq i<j\leq n}\abs{t_i-t_j}^{2\gamma}\,
\dup t_1 \cdots  \dup t_n \\
=\prod_{i,j=1}^n 
\frac{\Gamma(\alpha+(2n-i-j)+\la_i+\mu_j)}
{\Gamma(\alpha+(2n-i-j+1)+\la_i+\mu_j)}
\prod_{j=0}^{n-1} \frac{\Gamma((j+1)\gamma)\Gamma(1+(j+1)\gamma)}
{\Gamma(1+\gamma)}\\ \times
P_{\la}^{(1/\gamma)}(1^n) \,
P_{\mu}^{(1/\gamma)}(1^n).
\end{multline*}
This integral, which generalises the $\beta=\gamma$ case of
Kadell's integral \eqref{KadellInt} is originally due 
to Kadell \cite{Kadell93} and (for $\gamma=1$) L.K. Hua \cite{Hua79}.
Kadell's integral \eqref{KadellInt} also has an elliptic analogue,
which has the feature that the Dotsenko--Fateev integral  \eqref{DF} is
a special case
 \cite{Rains}. 

\medskip

There are other integrals is the literature referred to
as elliptic Selberg integrals, although they do not contain the
actual Selberg integral as a limiting case.
These integrals arise as solutions to the 
Knizhnik--Zamolodchikov--Bernard (KZB) heat equation
for $(2n+1)$-dimensional $\gsl_2$ modules
\[
2\pi\iup \kappa\, \frac{\partial u}{\partial\tau}=
\frac{\partial^2 u}{\partial\la^2}+n(n+1) \rho'(\la,\tau) u.
\]
Here $u=u(\la,\tau)$, $\rho(\la,\tau)=\vartheta'(\la,\tau)/
\vartheta(\la,\tau)$ with differentiation with respect to $\la$,
and $\vartheta(\la,\tau)=\theta_1(\pi\la,\tau)$ is a
theta function \cite{WW96}.

To describe the relevant solutions to the KZB equations let
$\Phi$ be the elliptic master function
\[
\Phi(t_1,\dots,t_n;\tau)=\prod_{i=1}^n E(t_i,\tau)^{-2n}
\prod_{1\leq i<j\leq n} E(t_i-t_j,\tau)^2,
\]
where $E(t,\tau)$ is the elliptic analogue of $t$
\[
E(t,\tau)=\frac{\vartheta(t,\tau)}{\vartheta'(0,\tau)}.
\]
The solutions considered by 
G.~Felder, L.~Stevens and Varchenko \cite{FSV03} are the
linear combinations
\[
u_{\kappa,m}(\la,\tau)=J_{\kappa,m}(\la,\tau)+
(-1)^{n+1} J_{\kappa,m}(-\la,\tau),
\]
where 
\[
J_{\kappa,m}(\la,\tau):=\Int_{0<t_n<\cdots<t_1<1}\!\!\!
\Phi^{1/\kappa}(t_1,\dots,t_n;\tau)
\theta_{\kappa,m}\Bigl(\la+\frac{2\abs{t}}{\kappa},\tau\Bigr)
\prod_{i=1}^n\sigma_{\la}(t_i,\tau)\,\dup t_1 \cdots \dup t_n.
\]
Whenever necessary this integral is understood in the sense of 
analytic continuation from the region where the exponents in 
$\Phi^{1/\kappa}$ have positive real part \cite{FSV03},
$\abs{t}=t_1+\cdots+t_n$, $\theta_{\kappa,m}(t,\tau)$ is a
theta function of degree $\kappa$ and characteristic $m$
\[
\theta_{\kappa,m}(\la,\tau)
=\sum_{j\in\Z+\frac{m}{2\kappa}} 
\eup^{2\pi\iup\kappa(\tau j+\la)j}
\]
and $\sigma_{\la}(t,\tau)=\theta(\la-t,\tau)/(\theta(\la,\tau) E(t,\tau))$.

In several instances Felder, Stevens and Varchenko found that the
``elliptic Selberg integrals'' $u_{\kappa,m}(\la,\tau)$ permit
closed form evaluations in terms of theta functions and ordinary
gamma functions. The simplest case of such an evaluation corresponds to
\begin{multline*}
u_{2n+2,n+1}(\la,\tau)=(2\pi)^{n/2} 
\eup^{-\pi\iup\frac{n(3n-1)}{4(n+1)}}
\eup^{\pi\iup\frac{n+1}{2}}
\theta(\la,\tau)^{n+1} \\ \times
S_n\biggl(\frac{n+2}{2(n+1)},-\frac{n}{n+1},\frac{1}{2(n+1)}\biggr)
\prod_{i=1}^n \Bigl(1-\eup^{2\pi\iup \frac{n+1+i}{2(n+1)}}\Bigr),
\end{multline*}
where $S_n$ is the Selberg integral \eqref{SelInt}.

\subsection*{The value distribution of 
\texorpdfstring{$\log\zeta(1/2+\iup t)$}{} on the critical line}

The final topic to be reviewed, following J.P.~Keating and N.C.~Snaith
\cite{KS00a}, is a link between the Selberg integral in its
trigonometric form \eqref{Morris}, and one of Selberg's theorems 
relating to the Riemann zeta function \cite{Selberg46}. 
The latter gives the value distribution
of $\log\zeta(1/2+\iup t)$ for large $t$, asserting that for 
any rectangle $B\in\Complex$
\begin{equation}\label{10.HD}
\lim_{T\to\infty} \frac{1}{T}
\biggl|\bigg\{ t: T \le t \le 2T, \: \:
\frac{\log \zeta (1/2 + \iup t)}{\sqrt{\frac{1}{2} \log \log T} }
\in B \bigg \} \biggr|
= \frac{1}{2} \iint_B \eup^{-\frac{1}{2}(x^2+y^2)} \, \dup x\, \dup y.
\end{equation}

To relate \eqref{10.HD} to \eqref{Morris}, first note that \eqref{MorrisEval}
can be used to evaluate
\begin{equation}\label{F1}
\Bigl\langle\; \prod_{i=1}^n \eup^{\frac{1}{2} \iup k \theta_i}
\abs{1+\eup^{\iup\theta_i}}^{\iup l} \Bigr \rangle,
\end{equation}
where the average is with respect to the eigenvalue PDF \eqref{Cp}
\cite{BF97}. Three further ingredients are then required:
the interpretation of \eqref{F1} as specifying a
distribution in random matrix theory, an hypothesis linking the
value distribution
of $\log \zeta(1/2+\iup t)$ to the value distribution of $\log \Lambda(z)$,
$\Lambda(z) := \prod_{i=1}^n \bigl(\exp(\iup\theta_i)-z\bigr)$ being the
characteristic polynomial for the random matrices, and the large $n$
form of \eqref{F1} deduced from \eqref{MorrisEval}.
Regarding the interpretation, note that
\begin{equation}\label{F2}
\Re \log \Lambda(-1) = \sum_{i=1}^n \log\abs{\eup^{\iup\theta_i}+1},
\qquad \Im \log \Lambda(-1) = \frac{1}{2} \sum_{i=1}^n \theta_i.
\end{equation}
It follows immediately that
\[
\int_{-\infty}^{\infty} \int_{-\infty}^{\infty}
\Bigl\langle \delta (s-\Re\log\Lambda(-1))
 \delta (t-\Im\log\Lambda(-1)) \Bigr\rangle 
\eup^{\iup l s} \eup^{\iup k t} \, \dup s\,\dup t
\]
is equal to \eqref{F1}. In other words, the characteristic function
for the joint distribution of the quantities \eqref{F2} is equal to
\eqref{F1}. 

The hypothesis of Keating and Snaith \cite{KS00a}, which extends the 
Montgomery--Odlyzko law linking the Riemann zeros to eigenvalues 
of large complex Hermitian random matrices (see e.g., \cite{KS99}), 
asserts that the value distribution of $\log \zeta(1/2+\iup t)$ 
for large $t$ will coincide with the value distribution of 
$\log \Lambda(z)$, $\abs{z}=1$, for $\Lambda(z)$ the characteristic
polynomial of matrices from the CUE ($\beta = 2$ case of \eqref{Cp}) 
for large $n$.
Further, the value of $n$ in the CUE is to be related to the value of $t$ in
$\zeta(1/2+\iup t)$ by $n=\log t$, which ensures that to leading order 
the density of eigenvalues and zeta function zeros is equal.

Thus the task at hand is to compute the large-$n$ limit of \eqref{F1} with
$k \mapsto k/(\frac{1}{\beta} \log n)^{1/2}$, 
$l \mapsto l/(\frac{1}{\beta} \log n)^{1/2}$,
which for $\beta = 2$ and with the identifications of the previous paragraph
corresponds to the scaling of $\log \zeta(1/2+\iup t)$ by 
$(\frac{1}{2} \log \log T)^{1/2}$ in \eqref{10.HD}. 
Using \eqref{MorrisEval} this limit has been computed in \cite{BF97} 
as being equal to $\exp(-(k^2 + l^2)/2)$. 
Hence the Selberg integral evaluation in its form
\eqref{MorrisEval} implies that the joint distribution of the scaled 
logarithm of the characteristic polynomial is equal to 
$\exp(-(s^2+t^2)/2)$, giving quantitive agreement between the 
hypothesis of Keating and Snaith and Selberg's theorem \eqref{10.HD}.

The value distribution of $\abs{\Lambda(z)}$ for $\abs{z}=1$ 
is also of relevance to zeta function theory \cite{KS00a}. 
The characteristic function of this quantity does not lead to a 
tractable integral. On the other hand, with
$p(s)$ a distribution supported on $s>0$, knowledge of the 
Mellin transform (complex moments)
\[
m(x) = \int_0^{\infty} s^{x-1} p(s) \, \dup s
\]
as a function in the complex plane gives, 
via the inverse Mellin transform,
\[
p(s) = \frac{1}{2 \pi \iup} \int_{c - \iup \infty}^{c + \iup \infty}
s^{-x} m(x) \, \dup x.
\]
With $p(s)$ the distribution of values of $\abs{\Lambda(-1)}^2$ for the CUE,
$m(x+1)$ is equal to \eqref{F1} with $k=0$, $\iup l = 2x$ and thus from
\eqref{MorrisEval} is given explicitly in terms of gamma functions
(for a discussion of computing the corresponding inverse Mellin transform, see
\cite{Rubinstein05}). It should also be remarked that the value distribution of
$\Lambda(\pm 1)$ for $\Lambda(z)$ the characteristic polynomial of a random
orthogonal or unitary symplectic matrix, chosen with Haar measure, is 
a special case of \eqref{BCS}, and thus similarly is an example of the
Selberg integral. Keating and Snaith \cite{KS00b} make use of this fact to
provide a quantitative link between the value distribution of families of
$L$-functions on the critical line and random matrix theory.


\begin{thebibliography}{99}

\bibitem{AIS07}
AMS--IMS--SIAM meeting
\textit{Interactions of Random Matrix Theory, Integrable Systems, 
and Stochastic Processes}, June 2007.

\bibitem{Anderson90}
G.W. Anderson,
\textit{The evaluation of Selberg sums},
C. R. Acad. Sci. Paris S\'er. I Math. \textbf{311} (1990), 469--472. 

\bibitem{Anderson91}
G.W. Anderson,
\textit{A short proof of Selberg's generalized beta formula},
Forum Math. \textbf{3} (1991), 415--417.

\bibitem{Andrews75}
G.E. Andrews,
\textit{Problems and prospects for basic hypergeometric functions},
in \textit{Theory and Application of Special Functions},
Math. Res. Center, Univ. Wisconsin, Publ. No. 35, 
Academic Press, New York, 1975; pp. 191--224.

\bibitem{AA81}
G. E. Andrews and R. Askey,
\textit{Another $q$-extension of the beta function},
Proc. Amer. Math. Soc. \textbf{81} (1981), 97--100. 

\bibitem{AAR99}
G. E. Andrews, R. Askey and R. Roy,
\textit{Special functions},
Encyclopedia of Mathematics and its Applications, Vol.~71,
Cambridge University Press, Cambridge, 1999.

\bibitem{Aomoto87}
K. Aomoto,
\textit{Jacobi polynomials associated with Selberg integrals},
SIAM J. Math. Anal. \textbf{18} (1987), 545--549.

\bibitem{Aomoto87b}
K. Aomoto,
\textit{The complex Selberg integral},
Quart. J. Math. Oxford  \textbf{38} (1987), 385--399.

\bibitem{Aomoto95}
K. Aomoto,
\textit{On a theta product formula for the symmetric A-type 
connection function},
Osaka J. Math. \textbf{32} (1995), 35--39.

\bibitem{Aomoto98}
K. Aomoto,
\textit{On elliptic product formulas for Jackson integrals 
associated with reduced root systems},
J. Algebraic Combin. \textbf{8} (1998), 115--126. 

\bibitem{Askey80} 
R. Askey,
\textit{Some basic hypergeometric extensions of integrals of
Selberg and Andrews},
SIAM J. Math. Anal. \textbf{11} (1980) 938--951.

\bibitem{Askey98}
R. Askey,
\textit{Letter to the SIAM minisymposium ``Problems and
solutions in special functions''}, in:
OP-SF NET 5.5 (Web resource), 1998.

\bibitem{AR89}
R. Askey and D. Richards, 
\textit{Selberg's second beta integral and an integral of Mehta}, in
\textit{Probability, Statistics, and Mathematics: 
Papers in Honor of Samuel Karlin}, T.W.~Anderson et al.~(eds), 
Academic Press, New York, 1989; pp.~27--39.

\bibitem{AW85}
R. Askey and J. Wilson,
\textit{Some basic hypergeometric orthogonal polynomials that 
generalize Jacobi polynomials},
Mem. Amer. Math. Soc. \textbf{54} (1985), Vol.~319.

\bibitem{BF97}
T.H. Baker and P.J. Forrester,
\textit{Finite-$N$ fluctuation formulas for random matrices},
J. Stat. Phys. \textbf{88} (1997), 1371--1386.

\bibitem{BF97b}
T.H. Baker and P.J. Forrester,
\textit{The Calogero--Sutherland model and generalized classical polynomials},
Comm. Math. Phys. \textbf{188} (1997), 175--216. 

\bibitem{BF98}
T.H. Baker and P.J. Forrester,
\textit{Nonsymmetric Jack polynomials and integral kernels},
Duke Math. J. \textbf{95} (1998), 1--50.  

\bibitem{Barnes04}
E.W. Barnes,
\textit{On the theory of the multiple gamma function},
Trans. Cambridge Phil. Soc. \textbf{19} (1904), 374--425.

\bibitem{BO93}
R.J. Beerends and E.M. Opdam,
\textit{Certain hypergeometric series related to the root system $BC$},
Trans. Amer. Math. Soc. \textbf{339} (1993), 581--609. 

\bibitem{Bergere98}
M.C. Berg\`ere,
\textit{Proof of Serban's conjecture},
J. Math. Phys. \textbf{39} (1998), 30--46.

\bibitem{Boas00}
R.P. Boas,
\textit{Comments on \cite{Polya00}}, in
\textit{George Polya, Collected Papers}, Vol.~1, 
R.P.~Boas (ed.), M.I.T.~Press, Cambridge, 1974; pp.~771--773.

\bibitem{Boas54}
R.P. Boas Jr.,
\textit{Book Reviews: New Journals}
Bull. Amer. Math. Soc. \textbf{60} (1954), 92--93.

\bibitem{Bombieri07}
E. Bombieri,
Private correspondence, 28 August 2007.

\bibitem{BO01}
A. Borodin and G. Olshanski,
\textit{$z$-Measures on partitions, Robinson--Schensted--Knuth
correspondence, and $\beta=2$ random matrix ensembles},
in \textit{Random matrix models and their applications},
P.M.~Bleher and A.R.~Its (eds), MSRI Publ. Vol.~40,
Cambridge University Press, Cambridge, 2001; pp.~71--94.

\bibitem{Calogero69}
F. Calogero,
\textit{Solution of a three-body problem in one dimension},
J. Math. Phys. \textbf{10} (1969), 2191--2196.

\bibitem{Carlson14}
F. Carlson, \textit{Sur une Classe de S\'eries de Taylor},
Ph.D. thesis, Uppsala Univ. 1914.

\bibitem{Cherednik91}
I.V. Cherednik,
\textit{A unification of the Knizhnik--Zamolodchikov and Dunkl
operators via affine Hecke algebras},
Inv. Math. \textbf{106} (1991), 411--432.

\bibitem{Cherednik95}
I.V. Cherednik,
\textit{Double affine Hecke algebras and Macdonald's conjectures},
Ann. of Math. \textbf{141} (1995), 191--216.

\bibitem{CR99}
P.B. Cohen and A. Regev,
\textit{Asymptotics of multinomial sums and identities 
between multi-integrals},
Israel J. Math. \textbf{112} (1999), 301--325. 

\bibitem{Constantine63}
A.G. Constantine,
\textit{Some noncentral distribution problems in multivariate analysis},
Ann. Math. Statist. \textbf{34} (1963), 1270--1285. 

\bibitem{CG06}
H. Coskun and R.A. Gustafson,
\textit{Well-poised Macdonald functions $W_{\lambda}$ and Jackson 
coefficients $\omega_{\lambda}$ on $BC_n$},
Contemp. Math. \textbf{417} (2006), 127--155.

\bibitem{Davis72}
A.W. Davis,
\textit{On the marginal distributions of the latent roots of the
multivariate beta matrix},
Ann. Math. Stat. \textbf{43} (1972), 1664--1669.

\bibitem{Debiard87}
A. Debiard,
\textit{Syst\`eme diff\'erentiel hyperg\'eom\'etrique et parties 
radiales des op\'erateurs invariants des espaces sym\'etriques 
de type $BC_p$},
Lecture Notes in Math. Vol.~1296, Springer, Berlin, 1987; pp. 42--124.

\bibitem{vanDiejen96}
J.F. van Diejen,
\textit{Self-dual Koornwinder--Macdonald polynomials},
Invent. Math. \textbf{126} (1996), 319--339. 

\bibitem{vanDiejen97}
J.F. van Diejen,
\textit{Confluent hypergeometric orthogonal polynomials related to 
the rational quantum Calogero system with harmonic confinement},
Comm. Math. Phys. \textbf{188} (1997), 467--497. 

\bibitem{DS01}
J.F. van Diejen and V.P. Spiridonov,
\textit{Elliptic Selberg integrals},
Internat. Math. Res. Notices (2001), 1083--1110.

\bibitem{Dixon05}
A.L. Dixon,
\textit{Generalizations of Legendre's formula $KE'-(K-E)K'=\frac{1}{2}\pi$},
Proc. London Math. Soc. \textbf{3} (1905), 206--224.

\bibitem{DF85}
V.S. Dotsenko and V.A. Fateev,
\textit{Four-point correlation functions and the operator algebra in
2D conformal invariant theories with central charge $C \le 1$},
Nucl. Phys. B \textbf{251} (1985), 691--734.

\bibitem{DE02}
I. Dumitriu and A. Edelman,
\textit{Matrix models for beta ensembles},
J. Math. Phys. \textbf{43} (2002), 5830--5847.

\bibitem{Dunkl89}
C.F. Dunkl,
\textit{Differential-difference operators associated to reflection
groups},
Trans. Amer. Math. Soc. \textbf{311}, 167--183.

\bibitem{DX01}
C.F. Dunkl and Y. Xu,
\textit{Orthogonal polynomials of several variables},
Encyclopedia of Mathematics and its Applications, Vol.~81,
Cambridge University Press, Cambridge, 2001.

\bibitem{Dyson62}
F.J. Dyson,
\textit{Statistical theory of energy levels of complex systems. I},
J. Math. Phys. \textbf{3} (1962), 140--156.

\bibitem{EFK03}
P.I. Etingof, I.B. Frenkel and A.A. Kirillov, Jr,
\textit{Lectures on Representation Theory and 
Knizhnik--Zamolodchikov Equations},
Mathematical Surveys and Monographs, Vol.~58, 
Amer. Math. Soc., Providence, RI, 2003.

\bibitem{Euler1730}
L. Euler,
\textit{De progressionibus transcendentibus seu quarum 
termini generales algebraice dari nequeunt},
Comm. Acad. Sci. Petropolitanae \textbf{5} (1730), 36--57.

\bibitem{Euler1769}
L. Euler,
\textit{Institutiones Calculi Integralis, II, Opera Omnia},
Ser. 1, Vol 12.

\bibitem{Evans81}
R.J. Evans,
\textit{Identities for products of Gauss sums over finite fields},
L'Enseignement Math. \textbf{27} (1981), 197--209.

\bibitem{Evans91}
R.J. Evans,
\textit{The evaluation of Selberg character sums},
L'Enseignement Math. \textbf{37} (1991), 235--248.

\bibitem{Evans92}
R.J. Evans,
\textit{Multidimensional $q$-beta integrals},
SIAM J. Math. Anal. \textbf{23} (1992), 758--765.

\bibitem{Evans94}
R.J. Evans,
\textit{Multidimensional beta and gamma integrals},
Contemp. Math. \textbf{166} (1994), 341--357.

\bibitem{Evans95}
R.J. Evans,
\textit{Selberg--Jack character sums of dimension $2$},
J. Number Theory \textbf{54} (1995), 1--11.

\bibitem{Evans07}
R.J. Evans,
Private correspondence, 31 August 2007.

\bibitem{FSV03}
G. Felder, L. Stevens and A. Varchenko,
\textit{Elliptic Selberg integrals and conformal blocks},
Math. Res. Lett. \textbf{10} (2003), 671--684.

\bibitem{Fo93}
P.J. Forrester,
\textit{Recurrence equations for the computation of correlations
in the $1/r^2$ quantum many body system},
J. Stat. Phys. \textbf{72} (1993), 39--50.

\bibitem{Fo02}
P.J. Forrester,
\textit{Log-Gases and Random Matrices},
\urlw{ms.unimelb.edu.au/~matpjf}

\bibitem{Fo06}
P.J. Forrester,
\textit{Beta random matrix ensembles},
to appear, Proceedings of the IMS (Singapore) programme on
Random Matrix Theory and its Applications to Statistics and Wireless
Communications.

\bibitem{Fo07}
P.J. Forrester,
\textit{A random matrix decimation procedure relating $\beta = 2/(r+1)$ to
$\beta = 2(r+1)$},
preprint.

\bibitem{FR02}
P.J. Forrester and E.M. Rains,
\textit{Interpretations of some parameter dependent generalizations of
classical matrix ensembles},
Probab. Theory Relat. Fields \textbf{131} (2005), 1--61.

\bibitem{FR07}
P.J. Forrester and E.M. Rains,
in preparation.

\bibitem{Garvan89}
F.G. Garvan,
\textit{Some Macdonald--Mehta integrals by brute force},
in \textit{$q$-Series and Partitions},
IMA Vol. Math. Appl. 18, Springer, New York, 1989; pp.~77--98.

\bibitem{Garvan}
F.G. Garvan,
Unpublished computer proof of the $k=1$ case of \eqref{ConG} for H$_4$.

\bibitem{Garvan90}
F.G. Garvan,
\textit{A proof of the Macdonald--Morris root system conjecture for $F\sb 4$},
SIAM J. Math. Anal. \textbf{21} (1990), 803--821. 

\bibitem{GG91}
F.G. Garvan and G. Gonnet,
\textit{Macdonald's constant term conjectures for exceptional 
root systems},
Bull. Amer. Math. Soc. (N.S.) \textbf{24} (1991), 343--347. 

\bibitem{GR04}
G. Gasper and M. Rahman,
\textit{Basic Hypergeometric Series},
Encyclopedia of Mathematics and its Applications, Vol.~35,
second edition,
Cambridge University Press, Cambridge, 2004.

\bibitem{Gelfond29}
A.O. Gelfond,
\textit{Sur un theor\`eme de M.G.~Polya},
Atti Reale Accad. Naz Lincei \textbf{10} (1929), 569--574.

\bibitem{Good70}
I.J. Good,
\textit{Short proof of a conjecture of Dyson},
J. Math. Phys. \textbf{11} (1970), 1884.

\bibitem{Gunson62}
J. Gunson,
\textit{Proof of a conjecture of Dyson in the statistical 
theory of energy levels},
J. Math. Phys. \textbf{3} (1962), 752--753.

\bibitem{Gustafson90}
R.A. Gustafson,
\textit{A generalization of the Selberg beta-integral},
Bull. Am. Math. Soc. \textbf{22} (1990), 97--105.

\bibitem{Gustafson94}
R.A. Gustafson, 
\textit{Some q-beta integrals on $SU(n)$ and $Sp(n)$ that generalize 
the Askey--Wilson and Nassrallah--Rahman integrals},
SIAM J. Math. Anal. \textbf{25} (1994), 441--449.

\bibitem{Habsieger86}
L. Habsieger,
\textit{La $q$-conjecture de Macdonald--Morris pour $G_2$},
C.R. Acad. Sc. Paris S\'er I \textbf{303} (1986), 211--214.

\bibitem{Habsieger88}
L. Habsieger,
\textit{Une $q$-int\'{e}grale de Selberg et Askey},
SIAM J. Math. Anal. \textbf{19} (1988), 1475--1489.

\bibitem{Heckman91}
G.J. Heckman,
\textit{An elementary approach to the hypergeometric shift 
operators of Opdam},
Invent. Math. \textbf{103} (1991), 341--350. 

\bibitem{Herz55}
C.S. Herz,
\textit{Bessel functions of matrix argument},
Ann. Math. \textbf{61} (1955), 474--523.

\bibitem{Hua79}
L.K. Hua,
\textit{Harmonic analysis of functions of several complex variables in
the classical domains},
Translations of Mathematical Monographs, Vol.~6,
AMS, Providence, RI, 1979.

\bibitem{IAS07}
Institute for Advanced Study,
institute for advanced study: press releases: Atle Selberg 1917--2007,
\urlw{ias.edu/newsroom/announcements/view/1186683853.html}

\bibitem{Ismail05}
M.E.H. Ismail,
\textit{Classical and Quantum Orthogonal Polynomials in One Variable},
Encyclopedia of Mathematics and its Applications, Vol.~98,
Cambridge University Press, Cambridge, 2005.

\bibitem{Ito97}
M. Ito,
\textit{On a theta product formula for Jackson integrals associated 
with root systems of rank two},
J. Math. Anal. Appl. \textbf{216} (1997), 122--163. 

\bibitem{Johansson00}
K. Johansson,
\textit{Shape fluctuations and random matrices},
Comm. Math. Phys. \textbf{209} (2000), 437--476.

\bibitem{Johansson02}
K. Johansson,
\textit{Non-intersecting paths, random tilings and random matrices},
Prob. Theory Rel. Fields \textbf{123} (2002), 225--280.

\bibitem{Kadell88a}
K.W.J. Kadell,
\textit{A proof of some $q$-analogues of Selberg's integral for $k=1$},
SIAM J. Math. Analysis \textbf{19} (1988), 944--968.

\bibitem{Kadell88}
K.W.J. Kadell,
\textit{A proof of Askey's conjectured $q$-analogue of Selberg's integral
and a conjecture of Morris},
SIAM J. Math. Analysis \textbf{19} (1988), 969--986.

\bibitem{Kadell93}
K.W.J. Kadell,
\textit{An integral for the product of two Selberg--Jack 
symmetric polynomials},
Compositio Math. \textbf{87} (1993), 5--43.

\bibitem{Kadell94}
K.W.J. Kadell,
\textit{A proof of the $q$-Macdonald--Morris conjecture for $BC_n$},
Mem. Amer. Math. Soc. \textbf{108} (1994), Vol.~516.

\bibitem{Kadell97}
K.W.J. Kadell,
\textit{The Selberg-Jack symmetric functions},
Adv. Math. \textbf{130} (1997), 33--102.

\bibitem{Kakei98}
S. Kakei,
\textit{Intertwining operators for a degenerate double affine 
Hecke algebra and multivariable orthogonal polynomials},
J. Math. Phys. \textbf{39} (1998), 4993--5006.

\bibitem{Kaneko93}
J. Kaneko,
\textit{Selberg integrals and hypergeometric functions associated
with Jack polynomials},
SIAM J. Math Anal. \textbf{24} (1993), 1086--1110.

\bibitem{Kaneko96}
J. Kaneko,
\textit{$q$-Selberg integrals and Macdonald polynomials},
Ann. Sci. \'Ecole Norm. Sup. (4) \textbf{29} (1996), 583--637.


\bibitem{KS53}
S. Karlin and L.S. Shapley,
\textit{Geometry of moment space},
Mem. Amer. Math. Soc. \textbf{1953} (1953), Vol.~12.

\bibitem{KS99}
N. Katz and P. Sarnak,
\textit{Zeroes of zeta functions and symmetry},
Bull. Amer. Math. Soc. \textbf{36} (1999), 1--26.

\bibitem{KS00a}
J.P. Keating and N.C. Snaith,
\textit{Random matrix theory and $\zeta(1/2 + it)$},
Comm. Math. Phys. \textbf{214} (2001), 57--89.

\bibitem{KS00b}
J.P. Keating and N.C. Snaith,
\textit{Random matrix theory and $L$-functions at $s=1/2$},
Comm. Math. Phys. \textbf{214} (2001), 91--110.

\bibitem{KLR03}
J.P. Keating, N. Linden and Z. Rudnick,
\textit{Random matrix theory, the exceptional Lie groups and
$L$-functions},
J. Phys. A \textbf{36} (2003), 2933--2944.

\bibitem{KN04}
R. Killip and I. Nenciu,
\textit{Matrix models for circular ensembles},
Int. Math. Res. Not. \textbf{50} (2004), 2665--2701.

\bibitem{Koornwinder92}
T.H. Koornwinder,
\textit{Askey--Wilson polynomials for root systems of type BC},
Contemp. Math. \textbf{138} (1992), 189--204.

\bibitem{Koranyi}
A. Kor\'anyi,
\textit{Hua-type integrals, hypergeometric functions and 
symmetric polynomials},
in \textit{International Symposium in Memory of Hua Loo Keng, Vol. II},
Springer, Berlin, 1991; pp. 169--180.

\bibitem{Lassalle91a}
M. Lassalle,
\textit{Polyn\^omes de Jacobi g\'en\'eralis\'es},
C. R. Acad. Sci. Paris S\'er. I Math. \textbf{312} (1991), 425--428.

\bibitem{Lassalle91b}
M. Lassalle,
\textit{Polyn\^omes de Laguerre g\'en\'eralis\'es},
C. R. Acad. Sci. Paris S\'er. I Math. \textbf{312} (1991), 725--728.

\bibitem{Lassalle91c}
M. Lassalle,
\textit{Polyn\^omes de Hermite g\'en\'eralis\'es},
C. R. Acad. Sci. Paris S\'er. I Math. \textbf{313} (1991), 579--582.

\bibitem{LT03}
J.-G. Luque and J.-Y. Thibon,
\textit{Hankel hyperdeterminants and Selberg integrals},
J. Phys. A \textbf{36} (2003), 5267--5292.

\bibitem{NR85}
B. Nassrallah and M. Rahman,
\textit{Projection formulas, a reproducing kernel and a 
generating function for $q$-Wilson polynomials},
SIAM J. Math. Anal. \textbf{16} (1985), 186--197.

\bibitem{Macdonald82}
I.G. Macdonald,
\textit{Some conjectures for root systems},
SIAM J. Math. Anal. \textbf{13} (1982), 988--1007.

\bibitem{Macdonald87}
I.G. Macdonald,
\textit{Commuting differential operators and zonal spherical functions}, 
Lecture Notes in Math. \textbf{1271} (1987), 189--200.

\bibitem{Macdonald95}
I.G. Macdonald,
\textit{Hall polynomials and symmetric functions}, 2nd ed.,
Oxford University Press, Oxford, 1995.

\bibitem{Macdonald03}
I.G. Macdonald,
\textit{Affine Hecke algebras and orthogonal polynomials},
Cambridge Tracts in Mathematics, Vol.~157, 
Cambridge University Press, Cambridge, 2003.

\bibitem{Macdonald}
I. G. Macdonald,
\textit{Hypergeometric functions}, unpublished manuscript.

\bibitem{Mehta67}
M.L. Mehta,
\textit{Random Matrices and the Statistical Theory of Energy Levels}, 
Academic Press, New York, 1967.

\bibitem{Mehta74}
M.L. Mehta,
\textit{Problem 74--6, Three multiple integrals},
SIAM Review \textbf{16} (1974), 256--257.

\bibitem{Mehta04}
M.L. Mehta,
\textit{Random Matrices}, 3rd ed., 
Pure and Applied Mathematics, Vol.~142,
Elsevier/Academic Press, Amsterdam, 2004.

\bibitem{MD63}
M.L. Mehta and F.J. Dyson,
\textit{Statistical theory of the energy levels of complex systems. V},
J. Math. Phys. \textbf{4} (1963), 713--719.

\bibitem{Mi93}
K. Mimachi,
\textit{Reducibility and irreducibility of the Gauss-Manin system associated
with a Selberg type integral}, \textbf{132} (1993), 43--62.

\bibitem{MY03}
K. Mimachi and M. Yoshida,
\textit{The reciprocity relation of the Selberg function},
in \textit{Proceedings of the International Conference on 
Special Functions and their Applications (Chennai, 2002)},
J. Comput. Appl. Math. \textbf{160} (2003), 209--215.

\bibitem{MY04}
K. Mimachi and M. Yoshida,
\textit{Intersection numbers of twisted cycles associated with the 
Selberg integral and an application to the conformal field theory},
Comm. Math. Phys. \textbf{250} (2004), 23--45.

\bibitem{MT05}
K. Mimachi and T. Takamuki,
\textit{A generalization of the beta integral arising from 
the Knizhnik--Zamolodchikov equation for the vector 
representations of types $B_n$, $C_n$ and $D_n$},
Kyushu J. Math. \textbf{59} (2005), 117--126. 

\bibitem{Morris82}
W.G. Morris,
\textit{Constant Term Identities for Finite and Affine Root Systems:
Conjectures and Theorems},
Ph.D. thesis, Univ. Wisconsin--Madison, 1982.

\bibitem{Muirhead70}
R. J. Muirhead,
\textit{Systems of partial differential equations for
hypergeometric functions of matrix argument},
Ann. Math. Statist. \textbf{41} (1970), 991--1001.

\bibitem{MV00}
E. Mukhin and A. Varchenko,
\textit{Remarks on critical points of phase functions and norms of
Bethe vectors},
Adv. Stud. Pure Math. \textbf{27} (2000), 239--246.

\bibitem{Okounkov98}
A. Okounkov,
(Shifted) Macdonald polynomials: $q$-integral representation and
combinatorial formula,
Compositio Math. \textbf{112} (1998), 147--182.

\bibitem{Okounkov98b}
A. Okounkov,
\textit{BC-type interpolation Macdonald polynomials and binomial formula 
for Koornwinder polynomials},
Transform. Groups \textbf{3} (1998) 181--207.

\bibitem{Okounkov01}
A. Okounkov,
\textit{Infinite wedge and random partitions},
Selecta Math., New Ser. \textbf{7} (2001), 57--81.

\bibitem{Opdam89}
E.M. Opdam,
\textit{Some applications of hypergeometric shift operators},
Invent. Math. \textbf{98} (1989), 1--18.

\bibitem{Opdam93}
E.M. Opdam,
\textit{Dunkl operators, Bessel functions and the discriminant of a
finite Coxeter group},
Compositio Math. \textbf{85} (1993), 333--373.

\bibitem{Opdam00}
E.M. Opdam,
\textit{Lecture notes on Dunkl operators for real and 
complex reflection groups},
MSJ Memoirs, Vol.~8, Math. Soc. Japan, Tokyo, 2000.

\bibitem{Polya00}
G. Polya,
\textit{Collected Papers, Vol.~1}, 
R.P.~Boas (ed.), M.I.T.~Press, Cambridge, 1974; pp.~1--16.

\bibitem{Rahman86}
M. Rahman,
\textit{An integral representation of a $_{10}\varphi_9$ 
and continuous bi-orthogonal $_{10}\varphi_9$ rational functions},
Canad. J. Math. \textbf{38} (1986), 605--618. 

\bibitem{Rains07}
E.M. Rains,
\textit{Transformations of elliptic hypergeometric integrals},
to appear in Ann. of Math.

\bibitem{Rains}
E.M. Rains,
Private communication.

\bibitem{Rains06}
E.M. Rains,
\textit{Limits of elliptic hypergeometric integrals},
to appear in Ramanujan J.

\bibitem{Regev81}
A. Regev,
\textit{Asymptotic values for degrees associated 
with strips of Young diagrams},
Adv. in Math. \textbf{41} (1981), 115--136.

\bibitem{Regev07}
A. Regev,
Private correspondence, 30 August 2007.

\bibitem{RZ02}
D. Richards and Q. Zheng,
\textit{Determinants of period matrices and an application to Selberg's
multidimensional beta integral},
Adv. in Appl. Math. \textbf{28} (2002), 602--633.

\bibitem{Rubinstein05}
M. Rubinstein,
\textit{Computational methods and experiments in analytic number theory},
\textit{Recent Perspectives in Random Matrix Theory and Number Theory},
F.~Mezzadri and N.C.~Snaith (eds), LMS Lecture Note Series \textbf{322},  
Cambridge University Press, Cambridge, 2005; pp.~425--506.

\bibitem{Ruijsenaars97}
S.N.M. Ruijsenaars,
\textit{First order analytic difference equations and integrable 
quantum systems},
J. Math. Phys. \textbf{38} (1997), 1069--1146.

\bibitem{SV91}
V.V. Schechtman and A.N. Varchenko, 
\textit{Arrangements of hyperplanes and Lie algebra homology},
Invent. Math. \textbf{106} (1991), 139--194. 

\bibitem{Se41}
A. Selberg,
\textit{\"Uber einen Satz von A. Gelfond}, 
Arch. Math. Naturvid. \textbf{44} (1941) 159--171.

\bibitem{Selberg44}
A. Selberg,
\textit{Bemerkninger om et multipelt integral},
Norsk. Mat. Tidsskr. \textbf{24} (1944), 71--78.

\bibitem{Selberg46}
A. Selberg,
\textit{Contributions to the theory of the Riemann zeta-function},
Arch. Math. OG. Naturv. B \textbf{48} (1946), 89--155.
(Reprinted with commentary in \cite{Selberg89}.) 

\bibitem{Selberg89}
A. Selberg,
\textit{Collected papers I},
Springer--Verlag, Heidelberg, 1989; p.~212.

\bibitem{Spiridonov01}
V.P. Spiridonov,
\textit{On the elliptic beta function},
Russ. Math. Surveys \textbf{56} (2001), 185--186.

\bibitem{Spiridonov01b}
V.P. Spiridonov,
\textit{Elliptic beta integrals and special functions of hypergeometric type},
in \textit{Integrable structures of exactly solvable two-dimensional models
of quantum field theory}, 
NATO Sci. Ser. II Math. Phys. Chem. Vol.~35, 
Kluwer Acad. Publ., Dordrecht; pp. 305--313.

\bibitem{Spiridonov07}
V.P. Spiridonov,
\textit{Short proofs of the elliptic beta integrals},
Ramanujan J. \textbf{13} (2007), 265--283. 

\bibitem{Spiridonov}
V.P. Spiridonov,
\textit{Elliptic hypergeometric functions},
A complement to \cite{AAR99}, written for its Russian edition.

\bibitem{Stanley89}
R.P. Stanley,
\textit{Some combinatorial properties of Jack symmetric functions},
Adv. Math. \textbf{77} (1989), 76--115. 

\bibitem{Stanley97}
R.P. Stanley,
\textit{Enumerative Combinatorics, Vol. I}, 
Cambridge University Press, New York/\-Cam\-bridge, 1997.

\bibitem{Stanley05}
R.P. Stanley,
\textit{Queue problems revisited},
Suomen Teht\"av\"aniekat \textbf{59} (2005), 193--203.

\bibitem{Stanley07}
R.P. Stanley,
Private correspondence, 15 September 2007.

\bibitem{Stembridge88}
J.R. Stembridge,
\textit{A short proof of Macdonald's conjecture for 
the root systems of type A},
Proc. Amer. Math. Soc. \textbf{102} (1988), 777--786.

\bibitem{Stokman97}
J.V. Stokman,
\textit{Multivariable big and little $q$-Jacobi polynomials},
SIAM J. Math. Anal. \textbf{28} (1997), 452--480. 

\bibitem{Stokman00}
J.V. Stokman,
\textit{On $BC$ type basic hypergeometric orthogonal polynomials},
Trans. Amer. Math. Soc. \textbf{352} (2000), 1527--1579. 


\bibitem{Sutherland71}
B. Sutherland, 
\textit{Exact results for a quantum many body problem in one-dimension},
Phys. Rev. A \textbf{4} (1971), 2019--2021.

\bibitem{Sutherland72}
B. Sutherland,
\textit{Exact results for a quantum many-body problem in one dimension: II},
Phys. Rev. A \textbf{5} (1972), 1372--1376.

\bibitem{TV00}
V. Tarasov and A. Varchenko, 
\textit{Difference equations compatible with trigonometric KZ differential 
equations},
Internat. Math. Res. Notices \textbf{2000}, 801--829.

\bibitem{TV03}
V. Tarasov and A. Varchenko, 
\textit{Selberg-type integrals associated with $\gsl_3$},
Lett. Math. Phys. \textbf{65} (2003), 173--185.

\bibitem{Varchenko89}
A. Varchenko,
\textit{The Euler beta-function, the Vandermonde determinant, the
Legendre equation, and critical values of linear functions on a
configuration of hyperplanes. I},
Math. USSR \textbf{35} (1990), 543--571.

\bibitem{Varchenko90}
A. Varchenko,
\textit{The Euler beta-function, the Vandermonde determinant, the
Legendre equation, and critical values of linear functions on a
configuration of hyperplanes. II},
Math. USSR \textbf{36} (1991), 155--167.

\bibitem{Varchenko03}
A. Varchenko, 
\textit{Special Functions, KZ Type Equations, and Representation Theory},
CBMS Regional Conference Series in Mathematics, Vol.~98,
Amer. Math. Soc., Providence, RI, 2003.

\bibitem{W05}
S.O. Warnaar,
\textit{$q$-Selberg integrals and Macdonald polynomials},
Ramanujan J. \textbf{10} (2005), 237--268.

\bibitem{W07}
S.O. Warnaar,
\textit{On the generalised Selberg integral of Richards and Zheng},
to appear in Adv. Appl. Math.

\bibitem{W07b}
S.O. Warnaar,
\textit{Bisymmetric functions, Macdonald polynomials and
$\gsl_3$ basic hypergeometric series},
to appear in Compositio Math.

\bibitem{W07c}
S.O. Warnaar,
\textit{A Selberg integral for the Lie algebra A$_n$},
arXiv:0708.1193.

\bibitem{Wilson62}
K. Wilson,
\textit{Proof of a conjecture of Dyson},
J. Math. Phys. \textbf{3} (1962), 1040--1043.

\bibitem{WW96}
E.T. Whittaker and G.N. Watson,
\textit{A course of modern analysis},
Reprint of the fourth (1927) edition,
Cambridge University Press, Cambridge, 1996.

\bibitem{Yan92}
Z. Yan,
\textit{A class of generalized hypergeometric functions in several variables},
Can J. Math. \textbf{44} (1992), 1317--1338.

\bibitem{ZB85}
D. Zeilberger and D.M. Bressoud,
\textit{A proof of Andrews' $q$-Dyson conjecture},
Discrete Math. \textbf{54} (1985), 201--224. 

\bibitem{Zeilberger87}
D. Zeilberger,
\textit{A proof of the $G_2$ case of Macdonald's root system-Dyson conjecture},
SIAM J. Math. Anal. \textbf{18} (1987), 880--883.

\bibitem{Zeilberger88}
D. Zeilberger,
\textit{A unified approach to Macdonald's root system conjecture},
SIAM J. Math. Anal. \textbf{19} (1988), 987--1013.

\bibitem{Zeilberger89}
D. Zeilberger,
\textit{A Stembridge--Stanton style proof of the Habsieger--Kadell 
$q$-Morris identity},
Discrete Math. \textbf{79} (1989), 313--322.

\end{thebibliography}
\end{document}